\documentstyle{article}

\newcommand{\eL}{\mbox{$\cal L$}}
\newcommand{\Pe}{\mbox{$\cal P$}}

\newcommand{\kon}{\wedge}
\newcommand{\dis}{\vee}
\newcommand{\str}{\rightarrow}
\newcommand{\rts}{\leftarrow}
\newcommand{\mj}{\mbox{\bf 1}}

\newcommand{\df}{\mbox{\scriptsize{\it df}}}
\newcommand{\HDS}{\vrule width0pt height2.3ex depth1.05ex\displaystyle}

\newcommand{\eM}{\mbox{$\cal M$}}

\newcommand{\mds}{\mbox{\bf DS}}

\newcommand{\cut}{\mbox{\it cut}}

\newcommand{\KR}{\mbox{$\kon^R$}}

\def\cirk{\,{\raisebox{.3ex}{\tiny $\circ$}}\,}
\def\ks{\mbox{\footnotesize$\;\xi\;$}}
\def\b#1#2{\stackrel{\raisebox{-2pt}{\mbox{\tiny $#1$}}}
{\raisebox{0pt}{$b$}}^{\raisebox{-7pt}{\scriptsize $#2$}}}

\def\c#1{\stackrel{\raisebox{-2pt}{\mbox{\tiny $\,#1$}}}
{\raisebox{0pt}{$c$}}}

\def\d#1#2{\stackrel{\raisebox{-2pt}{\mbox{\tiny $\,#1$}}}
{\raisebox{0pt}{$\delta$}}^{\raisebox{-7pt}{\scriptsize $#2$}}}

\def\s#1#2{\stackrel{\raisebox{-2pt}{\mbox{\tiny $#1$}}}
{\raisebox{0pt}{$\sigma$}}^{\raisebox{-7pt}{\scriptsize $#2$}}}

\def\f#1#2{{{\HDS #1}\over{\HDS #2}}}

\def\nav#1#2{\parshape=2 1em 31em 4em 28em \smallskip
\noindent{\makebox[3em][l]{#1}}{#2}\par\smallskip}

\def\prop#1#2{\vspace{2ex} \noindent{\sc #1.} {\it #2} \par \vspace{2ex}}
\def\dkz{\noindent{\sc Proof. }}
\def\qed{\hfill $\dashv$}

\def\pl{\!+\!}
\def\mn{\!-\!}

\def\e#1{\stackrel{\raisebox{-2pt}{\mbox{\tiny $\,#1$}}}
{\raisebox{0pt}{$e$}}}

\def\ep#1{\makebox(0,12){}\makebox(6,6)[b]{$\e{#1}$}'}

\def\den{=_{\mbox{\scriptsize{\it{dn}}}}}

\def\ro#1#2{\stackrel{\raisebox{-2pt}{\mbox{\tiny $\;#1$}}}
{\raisebox{0pt}{$\rho$}}^{\raisebox{-7pt}{\scriptsize $#2$}}}

\def\kst{\raisebox{1pt}{\mbox{\tiny$\xi$}}}

\newcommand{\PN}{\mbox{\bf PN}}
\newcommand{\PNN}{$\PN^{\neg}$}
\newcommand{\SMC}{\mbox{\bf SMC}}
\newcommand{\SA}{\mbox{\bf S}}

\newcommand{\SAps}{\mbox{\bf S}}

\newcommand{\SMCcons}{$\SMC^c$}
\newcommand{\SAcons}{$\SAps^c$}

\def\Del#1{\stackrel{\raisebox{-2pt}{\mbox{\tiny $#1$}}}
{\raisebox{0pt}{$\Delta$}}}
\def\Sig#1{\stackrel{\raisebox{-2pt}{\mbox{\tiny $#1$}}}
{\raisebox{0pt}{$\Sigma$}}}

\def\Dk{\Del{\kon}}
\def\Sd{\Sig{\dis}}
\def\Dkp{\Dk^{\raisebox{-7pt}{$'$}}}
\def\Sdp{\Sd^{\raisebox{-7pt}{$'$}}}

\def\Dd{\Del{\dis}}

\def\Sk{\Sig{\kon}}

\def\BE#1#2{\stackrel{\raisebox{-2pt}{\mbox{\tiny $\:#1$}}}
{\raisebox{0pt}{$B$}}^{\raisebox{-4.5pt}{\scriptsize $#2$}}}

\def\CE#1{\stackrel{\raisebox{-2pt}{\mbox{\tiny $\:#1$}}}
{\raisebox{0pt}{$C$}}}

\def\TE#1{\top^{#1}}

\def\NE#1{\bot^{\!#1}}

\def\NL{\neg^{\mbox{\scriptsize\rm L}}}

\def\NR{\neg^{\mbox{\scriptsize\rm R}}}

\def\ta#1{\stackrel{\raisebox{-2pt}{\mbox{\tiny $\,#1$}}}
{\raisebox{0pt}{$\tau$}}}

\def\tai#1{\stackrel{\raisebox{-2pt}{\mbox{\tiny $#1$}}}
{\raisebox{0pt}{$\tau$}}^{\raisebox{-7pt}{\scriptsize $-1$}}}

\def\Kon{\mbox{$\kon$}}

\def\Dis{\mbox{$\dis$}}

\def\koc{\begin{picture}(10,7)
\put(1.5,0){\line(0,1){7}} \put(8.5,0){\line(0,1){7}}
\put(1.5,0){\line(1,0){7}} \put(1.5,7){\line(1,0){7}}
\end{picture}}

\def\koci{\begin{picture}(7,5)
\put(1,0){\line(0,1){5}} \put(6,0){\line(0,1){5}}
\put(1,0){\line(1,0){5}} \put(1,5){\line(1,0){5}}
\end{picture}}

\begin{document}

\title{Coherence for Star-Autonomous Categories}
\author{{\sc Kosta Do\v sen} and {\sc Zoran Petri\' c}\\[0.5cm]
Mathematical Institute, SANU \\
Knez Mihailova 35, p.f. 367 \\
11001 Belgrade, Serbia \\
email: \{kosta, zpetric\}@mi.sanu.ac.yu}
\date{}
\maketitle

\begin{abstract}
\noindent This paper presents a coherence theorem for
star-autonomous categories exactly analogous to Kelly's and Mac
Lane's coherence theorem for symmetric monoidal closed categories.
The proof of this theorem is based on a categorial cut-elimination
result, which is presented in some detail.
\end{abstract}

\vspace{.3cm}


\vspace{.5ex}

\noindent {\it Keywords$\,$}: symmetric monoidal closed
categories, star-autonomous categories, cut elimination,
categorial coherence \vspace{1cm}

\baselineskip=1.2\baselineskip

\section{Introduction}

From the inception of proof nets in the late 1980s (see \cite{G87}
and \cite{DR89}), it could have been realized that they are
connected with the graphs one finds in Kelly's and Mac Lane's
coherence theorem for symmetric monoidal closed categories of
\cite{KML71}. The earliest explicit reference for that we know
about is \cite{Blute91} (see also \cite{Blute93}). It was also
soon suggested that the multiplicative fragment of classical
linear logic, which has an involutive negation that satisfies De
Morgan laws, is closely related to the notion of star-autonomous
category, which stems from \cite{Barr79} (see \cite{Laf88},
\cite{Seely89} and \cite{Barr91}).

Star-autonomous categories in the sense of \cite{Barr91} are
symmetric monoidal closed categories that have an object $\bot$
such that the canonical natural transformation from the identity
functor to the functor ${(\,\underline{\;\;\;}\str\bot)\str\bot}$
is a natural isomorphism (here
${\,\underline{\;\;\;}\str\underline{\;\;\;}\,}$ is the internal
hom-bifunctor). This notion is equivalent to the notion of
symmetric linearly (alias weakly) distributive category with
negation in the sense of \cite{CS97} (Section 4, Definition 4.3).
To establish the equivalence of the two notions is rather arduous,
as noted in \cite{CS97} (Theorem 4.5; a proof may be found in
\cite{DP05}, Chapter 3).

The aim of this paper is to present a coherence theorem for
symmetric linearly distributive categories with negation, which is
exactly analogous to Kelly's and Mac Lane's coherence theorem for
symmetric monoidal closed categories mentioned above. As Kelly's
and Mac Lane's proof of \cite{KML71}, the proof of our coherence
theorem is based on cut-elimination or similar results.  We will
not present all of them. Some of these results are in \cite{DP04},
and some in \cite{DP05} and \cite{DP05a}. We will present in some
detail only a cut-elimination theorem for symmetric linearly
distributive categories with negation freely generated by a set of
objects, on which our coherence theorem relies. This is a
cut-elimination theorem that asserts not only that for every
derivation we have a cut-free derivation of the same type, but the
original derivation and the cut-free derivation are moreover equal
as arrows in a category (which is not a preorder: not all arrows
of the same type are equal in this category).

As we indicated above, this paper is not self-contained. A more
detailed and more self-contained investigation of star-autonomous
categories and of their connection with the graphs of Kelly and
Mac Lane, and with the proof nets of classical linear logic, is in
the study \cite{DP05}.

Section 2, 3 and 5 of this paper introduce gradually the notion of
symmetric linearly distributive categories with negation freely
generated by a set of objects. Section 4 introduces a precise
notion of graph of the kind of Kelly and Mac Lane, and states the
previous coherence results on which we rely. Sections 6 and 7
contain the cut-elimination result, and Section 8 the coherence
result, which we have announced.

All the categories considered in this paper are small. We have no
need here for categories whose collections of objects or arrows
are bigger than sets.

\section{The category \mds}

The objects of the category \mds\ are the formulae of the
propositional language ${\eL_{\kon,\dis}}$, generated from a set
\Pe\ of propositional letters, which we call simply {\it letters},
with the binary connectives $\kon$ and $\dis$. We use
$p,q,r,\ldots\,$, sometimes with indices, for letters, and
$A,B,C,\ldots\,$, sometimes with indices, for formulae. As usual,
we omit the outermost parentheses of formulae and other
expressions later on.

To define the arrows of \mds, we define first inductively a set of
expressions called the {\it arrow terms} of \mds. Every arrow term
of \mds\ will have a {\it type}, which is an ordered pair of
formulae of ${\eL_{\kon,\dis}}$. We write ${f\!:A\vdash B}$ when
the arrow term $f$ is of type ${(A,B)}$. (We use the turnstile
$\vdash$ instead of the more usual $\str$, which we reserve for a
connective and a biendofunctor.) We use $f,g,h,\ldots\,$,
sometimes with indices, for arrow terms.

For all formulae $A$, $B$ and $C$ of ${\eL_{\kon,\dis}}$ the
following {\it primitive arrow terms}:

\begin{tabbing}
\centerline{$\mj_A\!: A\vdash A$,}
\\[1ex]
\mbox{\hspace{.2em}}\= $\b{\kon}{\str}_{A,B,C}\,$\= : \=
$A\kon(B\kon C)\vdash (A\kon B)\kon C$,\quad \=
$\b{\dis}{\str}_{A,B,C}\,$\= : \= $A\dis(B\dis C)\vdash (A\dis
B)\dis C$,\kill \> $\b{\kon}{\str}_{A,B,C}\,$\> : \> $A\kon(B\kon
C)\vdash (A\kon B)\kon C$,\> $\b{\dis}{\str}_{A,B,C}\,$\> : \>
$A\dis(B\dis C)\vdash (A\dis B)\dis C$,
\\*[1ex]
\> $\b{\kon}{\rts}_{A,B,C}\,$\> : \> $(A\kon B)\kon C\vdash A\kon
(B\kon C)$,\> $\b{\dis}{\rts}_{A,B,C}\,$\> : \> $(A\dis B)\dis
C\vdash A\dis (B\dis C)$,
\\[1ex]
\> \> $\c{\kon}_{A,B}\!\!$\' : \> $A\kon B\vdash B\kon A$,\> \>
$\c{\dis}_{A,B}\!\!$\' : \> $B\dis A\vdash A\dis B$,
\\[1ex]
\centerline{$d_{A,B,C}\!: A\kon(B\dis C)\vdash (A\kon B)\dis C$}
\end{tabbing}

\noindent are arrow terms of \mds. If ${g\!:A\vdash B}$ and
${f\!:B\vdash C}$ are arrow terms of \mds, then ${f\cirk
g\!:A\vdash C}$ is an arrow term of \mds; and if ${f\!:A\vdash D}$
and ${g\!:B\vdash E}$ are arrow terms of \mds, then ${f\ks
g\!:A\ks B\vdash D\ks E}$, for ${\!\ks\!\in\{\kon,\dis\}}$, is an
arrow term of \mds. This concludes the definition of the arrow
terms of \mds.

Next we define inductively the set of {\it equations} of \mds,
which are expressions of the form ${f=g}$, where $f$ and $g$ are
arrow terms of \mds\ of the same type. We stipulate first that all
instances of ${f=f}$ and of the following equations are equations
of \mds:

\begin{tabbing}
\mbox{\hspace{1em}}\= $({\mbox{{\it cat}~1}})$\quad\quad\= $f\cirk
\mj_A=\mj_B\cirk f=f\!:A\vdash B$,
\\*[1ex]
\> $({\mbox{{\it cat}~2}})$\> $h\cirk (g\cirk f)=(h\cirk g)\cirk
f$,
\\[1.5ex]
for $\!\ks\!\in\{\kon,\dis\}$,
\\*[1ex]
 \> $(\!\ks\, 1)$\> $\mj_A\ks\mj_B=\mj_{A\kst B}$,
\\*[1ex]
\> $(\!\ks\, 2)$\> $(g_1\cirk f_1)\ks(g_2\cirk f_2)=(g_1\ks
g_2)\cirk(f_1\ks f_2)$,
\\[1.5ex]
for $f\!:A\vdash D$, $g\!:B\vdash E$ and $h\!:C\vdash F$,
\\*[1ex]
\> $\mbox{($\b{\xi}{\str}$ {\it nat})}$\>  $((f\ks g)\ks
h)\cirk\!\b{\xi}{\str}_{A,B,C}\:=\:
\b{\xi}{\str}_{D,E,F}\!\cirk(f\ks (g\ks h))$,
\\[1ex]
\> $\mbox{($\c{\kon}$ {\it nat})}$\> $(g\kon
f)\cirk\!\c{\kon}_{A,B}\:=\:\c{\kon}_{D,E}\!\cirk(f\kon g)$,
\\[1ex]
\> $\mbox{($\c{\dis}$ {\it nat})}$\> $(g\dis
f)\cirk\!\c{\dis}_{B,A}\:=\:\c{\dis}_{E,D}\!\cirk(f\dis g)$,
\\[1ex]
\> $\mbox{($d$ {\it nat})}$\> $((f\kon g)\dis h)\cirk d_{A,B,C} =
d_{D,E,F}\cirk(f\kon (g\dis h))$,
\\[1.5ex]
\pushtabs $(\b{\xi}{}\b{\xi}{})$\quad\=
$\b{\xi}{\rts}_{A,B,C}\!\cirk\!\b{\xi}{\str}_{A,B,C}\;=\mj_{A\kst(B\kst
C)}$,\quad\quad\=
$\b{\xi}{\str}_{A,B,C}\!\cirk\!\b{\xi}{\rts}_{A,B,C}\;
=\mj_{(A\kst B)\kst C}$,
\\[1ex]
$(\b{\xi}{}\!5)$\> $\b{\xi}{\rts}_{A,B,C\kst D}\!\cirk\!
\b{\xi}{\rts}_{A\kst
B,C,D}\;=(\mj_A\:\ks\b{\xi}{\rts}_{B,C,D})\cirk\!\b{\xi}{\rts}_{A,B\kst
C,D}\!\cirk(\b{\xi}{\rts}_{A,B,C}\ks\:\mj_D)$,
\\[2ex]
$(\c{\kon}\c{\kon})$\>
$\c{\kon}_{B,A}\!\cirk\!\c{\kon}_{A,B}\;=\mj_{A\kon B}$,
\\*[1ex]
$(\c{\dis}\c{\dis})$\quad
$\c{\dis}_{A,B}\!\cirk\!\c{\dis}_{B,A}\;=\mj_{A\dis B}$,
\\[1.5ex]
$(\b{\kon}{}\c{\kon})$\> $(\mj_B\:\kon \c{\kon}_{C,A})\cirk\!
\b{\kon}{\rts}_{B,C,A}\!\cirk\!\c{\kon}_{A,B\kon
C}\!\cirk\!\b{\kon}{\rts}_{A,B,C}\!\cirk(\c{\kon}_{B,A}\kon\:\mj_C)=
\;\b{\kon}{\rts}_{B,A,C}$,
\\*[1ex]
$(\b{\dis}{}\c{\dis})$\> $(\mj_B\:\dis \c{\dis}_{A,C})\cirk\!
\b{\dis}{\rts}_{B,C,A}\!\cirk\!\c{\dis}_{B\dis C,A}
\!\cirk\!\b{\dis}{\rts}_{A,B,C}\!\cirk(\c{\dis}_{A,B}\dis\:\mj_C)=
\;\b{\dis}{\rts}_{B,A,C}$,
\\[1.5ex]
$(d \kon)$\> $(\b{\kon}{\rts}_{A,B,C}\dis\:\mj_D)\cirk d_{A\kon
B,C,D}=d_{A,B\kon C,D}\cirk(\mj_A\kon
d_{B,C,D})\cirk\!\b{\kon}{\rts}_{A,B,C\dis D}$,
\\[1ex]
$(d \dis)$\> $d_{D,C,B\dis
A}\cirk(\mj_D\:\kon\b{\dis}{\rts}_{C,B,A})=\;\b{\dis}{\rts}_{D\kon
C,B,A}\!\cirk(d_{D,C,B}\dis\mj_A)\cirk d_{D,C\dis B,A}$,
\\[1.5ex]
for $d^R_{C,B,A}=_{\df}\;\c{\dis}_{C,B\kon A} \!\cirk
(\c{\kon}_{A,B} \dis\:\mj_C)\cirk
d_{A,B,C}\cirk(\mj_A\:\kon\c{\dis}_{B,C})\cirk\!\c{\kon}_{C\dis
B,A}:$
\\*[.5ex]
\` $(C\dis B)\kon A\vdash C\dis(B\kon A)$,
\\[1ex]
$(d\!\b{\kon}{})$\> $d^R_{A\kon
B,C,D}\cirk(d_{A,B,C}\kon\mj_D)=d_{A,B,C\kon D}\cirk(\mj_A\kon
d^R_{B,C,D})\cirk\!\b{\kon}{\rts}_{A,B\dis C,D}$,
\\*[1ex]
$(d\!\b{\dis}{})$\>$(\mj_D\dis d_{C,B,A})\cirk d^R_{D,C,B\dis
A}=\;\b{\dis}{\rts}_{D,C\kon
B,A}\!\cirk(d^R_{D,C,B}\dis\mj_A)\cirk d_{D\dis C,B,A}$. \poptabs
\end{tabbing}

The set of equations of \mds\ is closed under symmetry and
transitivity of equality and under the rules

\[
(\mbox{\it cong~}\ks\!)\quad \f{f=f_1 \quad \quad \quad g=g_1}
{f\ks g=f_1\ks g_1}
\]

\noindent where ${\!\ks\!\in\{\cirk,\kon,\dis\}}$, and if
$\!\ks\!$ is $\cirk$, then ${f\cirk g}$ is defined (namely, $f$
and $g$ have appropriate, composable, types).

On the arrow terms of \mds\ we impose the equations of \mds. This
means that an arrow of \mds\ is an equivalence class of arrow
terms of \mds\ defined with respect to the smallest equivalence
relation such that the equations of \mds\ are satisfied (see
\cite{DP04}, Section 2.3, for details).

The equations $(\!\ks\, 1)$ and $(\!\ks\, 2)$ say that $\kon$ and
$\dis$ are biendofunctors (i.e.\ 2-endofunctors in the terminology
of \cite{DP04}, Section 2.4). Equations in the list above with
``\emph{nat}'' in their names, and analogous derivable equations,
will be called {\it naturality} equations. Such equations say that
$\b{\kon}{\str}$, $\b{\kon}{\rts}$, $\c{\kon}$, etc.\ are natural
transformations.

The equations $(d \kon)$, $(d \dis)$, ${(d\!\b{\kon}{})}$ and
${(d\!\b{\dis}{})}$ stem from \cite{CS97} (Section 2.1; see
\cite{CS92}, Section 2.1, for an announcement). The equation
${(d\!\b{\dis}{})}$ of \cite{DP04} (Section 7.2) amounts with
${(\b{\dis}{}\b{\dis}{})}$ to the present one.

\section{The category \PNN}

The category \PNN\ is defined as \mds\ save that we make the
following changes and additions. Instead of ${\eL_{\kon,\dis}}$,
we have the propositional language ${\eL_{\neg,\kon,\dis}}$, which
has in addition to what we have for ${\eL_{\kon,\dis}}$ the unary
connective $\neg$.

To define the arrow terms of \PNN, in the inductive definition we
had for the arrow terms of \mds\ we assume in addition that for
all formulae $A$ and $B$ of ${\eL_{\neg,\kon,\dis}}$ the following
{\it primitive arrow terms}:

\begin{tabbing}
\centerline{$\Dk_{B,A}:A\vdash A\kon(\neg B\dis B)$,}
\\*[1ex]
\centerline{$\Sd_{B,A}:(B\kon\neg B)\dis A\vdash A$,}
\end{tabbing}

\noindent are arrow terms of \PNN.

To define the arrows of \PNN, we assume in the inductive
definition we had for the equations of \mds\ the following
additional equations:

\begin{tabbing}
\mbox{\hspace{2em}}\= $\mbox{($\Dk$ {\it nat})}$\quad\=
$(f\kon\mj_{\neg B\dis B})\cirk\!\Dk_{B,A}\;=\;\Dk_{B,D}\!\cirk
f$,
\\*[1ex]
\> $\mbox{($\Sd$ {\it nat})}$\> $f\cirk\!
\Sd_{B,A}\;=\;\Sd_{B,D}\!\cirk(\mj_{B\kon \neg B}\dis f)$,
\\[2ex]
\> $\mbox{($\b{\kon}{}\Dk$)}$\> $\b{\kon}{\rts}_{A,B,\neg C\dis
C}\!\cirk\!\Dk_{C,A\kon B}\;=\mj_A\:\kon\Dk_{C,B}$,
\\*[1ex]
\> $\mbox{($\b{\dis}{}\Sd$)}$\> $\Sd_{C,B\dis
A}\!\cirk\!\b{\dis}{\rts}_{C\kon\neg
C,B,A}\;=\;\Sd_{C,B}\dis\:\mj_A$,
\\[2ex]
\quad for $\Sk_{B,A}\;=_{\df}\;\:\c{\kon}_{A,\neg B\dis
B}\!\cirk\!\Dk_{B,A}\: : A\vdash (\neg B\dis B)\kon A$,
\\*[1.5ex]
\> $\mbox{($d\!\Sk$)}$\> $d_{\neg A\dis A,B,C}\cirk\!\Sk_{A,B\dis
C}\;=\;\Sk_{A,B}\dis\:\mj_C$,
\\[2ex]
\quad for $\Dd_{B,A}\;
=_{\df}\;\:\Sd_{B,A}\!\cirk\!\c{\dis}_{B\kon\neg B,A}\:
:A\dis(B\kon\neg B)\vdash A$,
\\*[1.5ex]
\> $\mbox{($d\!\Dd$)}$\> $\Dd_{A,C\kon B}\!\cirk d_{C,B,A\kon\neg
A}\;=\;\mj_C\:\kon\Dd_{A,B}$,
\\[2ex]
\> $\mbox{($\Sd\Dk$)}$\> $\Sd_{A,A}\!\cirk d_{A,\neg A,A}
\cirk\!\Dk_{A,A}\;=\;\mj_A$,
\\*[2ex]
\quad for \=$\Dkp_{B,A}\;=_{\df}\;(\mj_A\:\kon\c{\dis}_{B,\neg
B})\cirk\!\Dk_{B,A}\: :A\vdash A\kon(B\dis\neg B)$ and
\\*[1ex]
\>$\Sdp_{B,A}\;=_{\df}\;\:\Sd_{B,A}\!\cirk(\c{\kon}_{\neg
B,B}\dis\:\mj_A):(\neg B\kon B)\dis A\vdash A$,
\\*[2ex]
\mbox{\hspace{2em}}$\mbox{($\Sdp\Dkp$)}$\>\> $\Sdp_{A,\neg
A}\!\cirk d_{\neg A,A,\neg A}\cirk\!\Dkp_{A,\neg A}\;=\mj_{\neg
A}$.
\end{tabbing}

The naturality equations \mbox{($\Dk$ {\it nat})} and \mbox{($\Sd$
{\it nat})} say that $\Dk$ and $\Sd$ are natural transformations
in the second index. We have analogous naturality equations for
$\Sk$, $\Dd$, $\Dkp$ and $\Sdp$.

The arrow ${\Dk_{B,A}:A\vdash A\kon(\neg B\dis B)}$ is analogous
to the arrow of type $A\vdash A\kon\top$ that one finds in
monoidal categories. However, $\Dk_{B,A}$ does not have an inverse
in \PNN. The equation $\mbox{($\b{\kon}{}\Dk$)}$ is analogous to
an equation that holds in monoidal categories (see \cite{ML71},
Section VII.1, \cite{DP04}, Section 4.6, and \S 5 below).

A {\it proof-net} category is a category with two biendofunctors
$\kon$ and $\dis$, a unary operation $\neg$ on objects, and the
natural transformations $\b{\kon}{\str}$, $\b{\kon}{\rts}$,
$\b{\dis}{\str}$, $\b{\dis}{\rts}$, $\c{\kon}$, $\c{\dis}$, $d$,
$\Dk$ and $\Sd$ that satisfy the equations $(\b{\xi}{}\!5)$,
$(\b{\xi}{}\b{\xi}{})$, $\ldots$ , $\mbox{($\Sdp\Dkp$)}$ of \PNN.

It is clear how to define the notion of proof-net functor between
proof-net categories, which preserves the proof-net structure of a
category strictly (i.e.\ ``on the nose''; cf.\ \cite{DP04},
Section 2.8). The functor $G$ from \PNN\ to \emph{Br} defined in
the next section is a proof-net functor in this sense. The other
functors $G$ mentioned later in the paper also preserve each a
certain categorial structure ``on the nose''.

The category \PNN\ is, up to isomorphism, the free proof-net
category generated by the set of letters \Pe, thought of as a
discrete category.

\section{The category \emph{Br}}

We are now going to introduce a category called \emph{Br}. This
category serves to formulate a coherence result for proof-net
categories, which says that there is a faithful functor from \PNN\
to \emph{Br}. The name of the category \emph{Br} comes from
``Brauerian''. The arrows of this category correspond to graphs,
or diagrams, that were introduced in \cite{B37} in connection with
Brauer algebras. Analogous graphs were investigated in
\cite{EK66}, and in \cite{KML71} Kelly and Mac Lane relied on them
to prove their coherence result for symmetric monoidal closed
categories.

Let \eM\ be a set whose subsets are denoted by $X$, $Y$,
$Z$,~$\ldots$ For ${i\in\{s,t\}}$ (where $s$ stands for ``source''
and $t$ for ``target''), let $\eM^i$ be a set in one-to-one
correspondence with \eM, and let ${i\!:\eM\str \eM^i}$ be a
bijection. Let $X^i$ be the subset of $\eM^i$ that is the image of
the subset $X$ of \eM\ under $i$. If ${u\in\eM}$, then we use
$u_i$ as an abbreviation for $i(u)$. We assume also that \eM,
$\eM^s$ and $\eM^t$ are mutually disjoint.

For ${X,Y\subseteq \eM}$, let a {\it split relation} of \eM\ be a
triple ${\langle R,X,Y\rangle}$ such that ${R\subseteq(X^s\cup
Y^t)^2}$. The set ${X^s\cup Y^t}$ may be conceived as the disjoint
union of $X$ and $Y$. We denote a split relation ${\langle
R,X,Y\rangle}$ more suggestively by ${R\!:X\vdash Y}$.

A split relation ${R\!:X\vdash Y}$ is a {\it split equivalence}
when $R$ is an equivalence relation. We denote by part($R$) the
partition of ${X_s\cup Y_t}$ corresponding to the split
equivalence ${R\!:X\vdash Y}$.

A split equivalence ${R\!:X\vdash Y}$ is {\it Brauerian} when
every member of part($R$) is a two-element set. For ${R\!:X\vdash
Y}$ a Brauerian split equivalence, every member of part($R$) is
either of the form ${\{u_s,v_t\}}$, in which case it is called a
{\it transversal}, or of the form ${\{u_s,v_s\}}$, in which case
it is called a {\it cup}, or, finally, of the form
${\{u_t,v_t\}}$, in which case it is called a {\it cap}.

For ${X,Y,Z\in\eM}$, we want to define the composition ${P\ast
R\!:X\vdash Z}$ of the split relations ${R\!:X\vdash Y}$ and
${P\!:Y\vdash Z}$ of \eM. For that we need some auxiliary notions.

For ${X,Y\subseteq \eM}$, let the function ${\varphi^s\!: X\cup
Y^t\str X^s\cup Y^t}$ be defined by

\[
\varphi^s(u)=\left\{
\begin{array}{ll}
u_s & \mbox{\rm{if }} u\in X
\\
u       & \mbox{\rm{if }} u\in Y^t,
\end{array}
\right.
\]

\noindent and let the function ${\varphi^t\!:X^s\cup Y\str X^s\cup
Y^t}$ be defined by

\[
\varphi^t(u)=\left\{
\begin{array}{ll}
u            & \mbox{\rm{if }} u\in X^s
\\
u_t       & \mbox{\rm{if }} u\in Y.
\end{array}
\right.
\]

For a split relation ${R\!:X\vdash Y}$, let the two relations
${R^{-s}\subseteq(X\cup Y^t)^2}$ and ${R^{-t}\subseteq(X^s\cup
Y)^2}$ be defined by

\[
(u,v)\in R^{-i}\quad\mbox{\rm{iff}}\quad
(\varphi^i(u),\varphi^i(v))\in R
\]

\noindent for ${i\in\{s,t\}}$. Finally, for an arbitrary binary
relation $R$, let $\mbox{\rm Tr}(R)$ be the transitive closure of
$R$.

Then we define ${P\ast R}$ by

\[
P\ast R=_{\df}\mbox{\rm Tr}(R^{-t}\cup P^{-s})\cap(X^s\cup Z^t)^2.
\]

\noindent It is easy to conclude that ${P\ast R\!:X\vdash Z}$ is a
split relation of \eM, and that if ${R\!:X\vdash Y}$ and
${P\!:Y\vdash Z}$ are (Brauerian) split equivalences, then ${P\ast
R}$ is a (Brauerian) split equivalence.

We now define the category \emph{Br}. The set of objects of
\emph{Br} is {\boldmath$N$}, the set of finite ordinals. The
arrows of \emph{Br} are the Brauerian split equivalences
${R\!:m\vdash n}$ of {\boldmath$N$}. The identity arrow
${\mj_n\!:n\vdash n}$ of \emph{Br} is the Brauerian split
equivalence such that

\[
\mbox{\rm{part}}(\mj_n)=\{\{m_s,m_t\}\mid m<n\}.
\]

\noindent Composition in \emph{Br} is the operation $\ast$ defined
above.

That \emph{Br} is indeed a category (i.e.\ that $\ast$ is
associative and that $\mj_n$ is an identity arrow) is proved in
\cite{DP03a} and \cite{DP03b}. This proof is obtained via an
isomorphic representation of \emph{Br} in the category \emph{Rel},
whose objects are the finite ordinals and whose arrows are all the
relations between these objects. Composition in \emph{Rel} is the
ordinary composition of relations. A direct formal proof would be
more involved, though what we have to prove is rather clear if we
represent Brauerian split equivalences geometrically (as this is
done in \cite{B37} and \cite{EK66}).

For example, for ${R\subseteq(3^s\cup 9^t)^2}$ and
${P\subseteq(9^s\cup 1^t)^2}$ such that

\begin{tabbing}
\hspace{2em}\=$\mbox{\rm{part}}(R)\;$\=$=\{\{0_s,0_t\},\{1_s,3_t\},\{2_s,6_t\}\}\cup
\{\{n_t,(n\pl 1)_t\}\mid n\in\{1,4,7\}\}$,
\\*[1.5ex]
\>$\mbox{\rm{part}}(P)$\>$=\{\{2_s,0_t\}\}\cup\{\{n_s,(n\pl
1)_s\}\mid n\in\{0,3,5,7\}\}$,
\end{tabbing}

\noindent the composition ${P\ast R\subseteq(3^s\cup 1^t)^2}$, for
which we have

\[
\mbox{\rm{part}}(P\ast R)={\{\{0_s,0_t\}, \{1_s,2_s\}\}},
\]

\noindent is obtained from the following diagram:

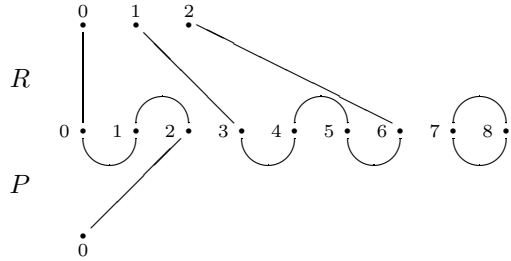
\begin{center}
\begin{picture}(160,115)

\put(3,23){\line(1,1){34}} \put(0,63){\line(0,1){34}}
\put(57,63){\line(-1,1){34}} \put(117,63){\line(-2,1){74}}

\put(0,20){\circle*{2}} \put(0,60){\circle*{2}}
\put(20,60){\circle*{2}} \put(40,60){\circle*{2}}
\put(60,60){\circle*{2}} \put(80,60){\circle*{2}}
\put(100,60){\circle*{2}} \put(120,60){\circle*{2}}
\put(140,60){\circle*{2}} \put(160,60){\circle*{2}}
\put(0,100){\circle*{2}} \put(20,100){\circle*{2}}
\put(40,100){\circle*{2}}

\put(10,57){\oval(20,20)[b]} \put(30,63){\oval(20,20)[t]}
\put(70,57){\oval(20,20)[b]} \put(110,57){\oval(20,20)[b]}
\put(90,63){\oval(20,20)[t]} \put(150,57){\oval(20,20)[b]}
\put(150,63){\oval(20,20)[t]}

\put(0,17){\makebox(0,0)[t]{\scriptsize$0$}}
\put(-5,60){\makebox(0,0)[r]{\scriptsize$0$}}
\put(15,60){\makebox(0,0)[r]{\scriptsize$1$}}
\put(35,60){\makebox(0,0)[r]{\scriptsize$2$}}
\put(55,60){\makebox(0,0)[r]{\scriptsize$3$}}
\put(75,60){\makebox(0,0)[r]{\scriptsize$4$}}
\put(95,60){\makebox(0,0)[r]{\scriptsize$5$}}
\put(115,60){\makebox(0,0)[r]{\scriptsize$6$}}
\put(135,60){\makebox(0,0)[r]{\scriptsize$7$}}
\put(155,60){\makebox(0,0)[r]{\scriptsize$8$}}
\put(0,103){\makebox(0,0)[b]{\scriptsize$0$}}
\put(20,103){\makebox(0,0)[b]{\scriptsize$1$}}
\put(40,103){\makebox(0,0)[b]{\scriptsize$2$}}

\put(-20,80){\makebox(0,0)[r]{$R$}}
\put(-20,40){\makebox(0,0)[r]{$P$}}

\end{picture}
\end{center}

\vspace{-2ex}

Every bijection $f$ from $X^s$ to $Y^t$ corresponds to a Brauerian
split equivalence ${R\!:X\vdash Y}$ such that the members of
part($R$) are of the form ${\{u,f(u)\}}$. The composition of such
Brauerian split equivalences, which correspond to bijections, is
then a simple matter: it amounts to composition of these
bijections. If in \emph{Br} we keep as arrows only such Brauerian
split equivalences, then we obtain a subcategory of \emph{Br}
isomorphic to the category \emph{Bij} whose objects are again the
finite ordinals and whose arrows are the bijections between these
objects. The category \emph{Bij} is a subcategory of the category
\emph{Rel} (which played an important role in \cite{DP04}), whose
objects are the finite ordinals and whose arrows are all the
relations between these objects. Composition in \emph{Bij} and
\emph{Rel} is the ordinary composition of relations. The category
\emph{Rel} (which played an important role in \cite{DP04}) is
isomorphic to a subcategory of the category whose arrows are split
relations of finite ordinals, of whom \emph{Br} is also a
subcategory.

We define a functor $G$ from \PNN\ to \emph{Br} in the following
way. On objects, we stipulate that $GA$ is the number of
occurrences of letters in $A$. On arrows, we have first that
${G\alpha}$ is an identity arrow of \emph{Br} for $\alpha$ being
$\mj_A$, $\b{\xi}{\str}_{A,B,C}$, $\b{\xi}{\rts}_{A,B,C}$ and
$d_{A,B,C}$, where ${\!\ks\!\in\{\kon,\dis\}}$.

Next, for ${i,j\in\{s,t\}}$, we have that ${\{m_i,n_j\}}$ belongs
to part(${G\!\c{\kon}_{A,B}}$) iff ${\{n_i,m_j\}}$ belongs to
part(${G\!\c{\dis}_{A,B}}$), iff $i$ is $s$ and $j$ is $t$, while
${m,n<GA\pl GB}$ and

\vspace{-1ex}

\[
(m\mn n\mn GA)(m\mn n\pl GB)=0.
\]

\noindent In the following example, we have ${G(p\dis
q)=2=\{0,1\}}$ and $G((q\dis\neg r)\dis q)$$=3=\{0,1,2\}$, and we
have the diagrams

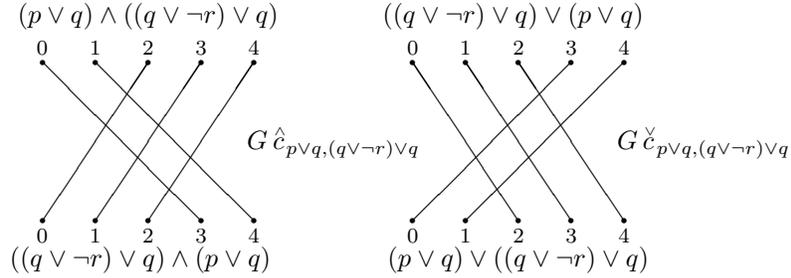
\vspace{2ex}

\begin{center}
\begin{picture}(300,95)
\put(10,10){\line(2,3){40}} \put(30,10){\line(2,3){40}}
\put(50,10){\line(2,3){40}} \put(70,10){\line(-1,1){60}}
\put(90,10){\line(-1,1){60}}

\put(150,10){\line(1,1){60}} \put(170,10){\line(1,1){60}}
\put(190,10){\line(-2,3){40}} \put(210,10){\line(-2,3){40}}
\put(230,10){\line(-2,3){40}}

\put(10,10){\circle*{2}} \put(30,10){\circle*{2}}
\put(50,10){\circle*{2}} \put(70,10){\circle*{2}}
\put(90,10){\circle*{2}} \put(150,10){\circle*{2}}
\put(170,10){\circle*{2}} \put(190,10){\circle*{2}}
\put(210,10){\circle*{2}} \put(230,10){\circle*{2}}

\put(10,70){\circle*{2}} \put(30,70){\circle*{2}}
\put(50,70){\circle*{2}} \put(70,70){\circle*{2}}
\put(90,70){\circle*{2}} \put(150,70){\circle*{2}}
\put(170,70){\circle*{2}} \put(190,70){\circle*{2}}
\put(210,70){\circle*{2}} \put(230,70){\circle*{2}}

\put(10,7){\makebox(0,0)[t]{\footnotesize 0}}
\put(30,7){\makebox(0,0)[t]{\footnotesize 1}}
\put(50,7){\makebox(0,0)[t]{\footnotesize 2}}
\put(70,7){\makebox(0,0)[t]{\footnotesize 3}}
\put(90,7){\makebox(0,0)[t]{\footnotesize 4}}
\put(150,7){\makebox(0,0)[t]{\footnotesize 0}}
\put(170,7){\makebox(0,0)[t]{\footnotesize 1}}
\put(190,7){\makebox(0,0)[t]{\footnotesize 2}}
\put(210,7){\makebox(0,0)[t]{\footnotesize 3}}
\put(230,7){\makebox(0,0)[t]{\footnotesize 4}}

\put(10,73){\makebox(0,0)[b]{\footnotesize 0}}
\put(30,73){\makebox(0,0)[b]{\footnotesize 1}}
\put(50,73){\makebox(0,0)[b]{\footnotesize 2}}
\put(70,73){\makebox(0,0)[b]{\footnotesize 3}}
\put(90,73){\makebox(0,0)[b]{\footnotesize 4}}
\put(150,73){\makebox(0,0)[b]{\footnotesize 0}}
\put(170,73){\makebox(0,0)[b]{\footnotesize 1}}
\put(190,73){\makebox(0,0)[b]{\footnotesize 2}}
\put(210,73){\makebox(0,0)[b]{\footnotesize 3}}
\put(230,73){\makebox(0,0)[b]{\footnotesize 4}}

\put(120,40){\makebox(0,0){$G\!\c{\kon}_{p\dis q,(q\dis\neg r)\dis
q}$}}

\put(260,40){\makebox(0,0){$G\!\c{\dis}_{p\dis q,(q\dis\neg r)\dis
q}$}}

\put(50,85){\makebox(0,0)[b]{$(p\dis q)\kon((q\dis\neg r)\dis
q)$}}

\put(47,0){\makebox(0,0)[t]{$((q\dis\neg r)\dis q)\kon(p\dis q)$}}

\put(188,85){\makebox(0,0)[b]{$((q\dis\neg r)\dis q)\dis(p\dis
q)$}}

\put(190,0){\makebox(0,0)[t]{$(p\dis q)\dis((q\dis\neg r)\dis
q)$}}

\end{picture}
\end{center}

\vspace{3ex}

We have that ${\{m_i,n_j\}}$ belongs to part(${G\!\Dk_{B,A}}$) iff
either

\nav{}{$i$ is $s$ and $j$ is $t$, while ${m,n<GA}$ and ${m=n}$,
or}

\vspace{-1ex}

\nav{}{$i$ and $j$ are both $t$, while ${m,n\in\{GA,\ldots,GA\pl
2GB\mn 1\}}$ and ${|m\mn n|=GB}$.}

\noindent In the following example, for $A$ being ${(q\dis\neg
r)\dis q}$ and $B$ being ${p\dis q}$, we have

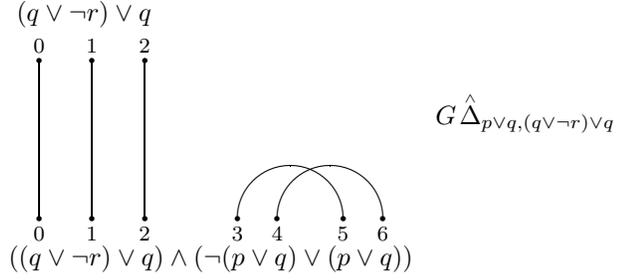
\vspace{2ex}

\begin{center}
\begin{picture}(180,107)
\put(10,20){\line(0,1){60}} \put(30,20){\line(0,1){60}}
\put(50,20){\line(0,1){60}}

\put(10,20){\circle*{2}} \put(30,20){\circle*{2}}
\put(50,20){\circle*{2}} \put(85,20){\circle*{2}}
\put(100,20){\circle*{2}} \put(125,20){\circle*{2}}
\put(140,20){\circle*{2}}

\put(10,80){\circle*{2}} \put(30,80){\circle*{2}}
\put(50,80){\circle*{2}}

\put(10,17){\makebox(0,0)[t]{\footnotesize 0}}
\put(30,17){\makebox(0,0)[t]{\footnotesize 1}}
\put(50,17){\makebox(0,0)[t]{\footnotesize 2}}
\put(85,17){\makebox(0,0)[t]{\footnotesize 3}}
\put(100,17){\makebox(0,0)[t]{\footnotesize 4}}
\put(125,17){\makebox(0,0)[t]{\footnotesize 5}}
\put(140,17){\makebox(0,0)[t]{\footnotesize 6}}

\put(10,83){\makebox(0,0)[b]{\footnotesize 0}}
\put(30,83){\makebox(0,0)[b]{\footnotesize 1}}
\put(50,83){\makebox(0,0)[b]{\footnotesize 2}}

\put(105,20){\oval(40,40)[t]} \put(120,20){\oval(40,40)[t]}

\put(160,60){\makebox(0,0)[l]{$G\!\Dk_{p\dis q,(q\dis\neg r)\dis
q}$}}

\put(75,10){\makebox(0,0)[t]{$((q\dis\neg r)\dis q)\kon(\neg(p\dis
q) \dis(p\dis q))$}}

\put(27,95){\makebox(0,0)[b]{$(q\dis\neg r)\dis q$}}

\end{picture}
\end{center}

\vspace{2ex}

We have that ${\{m_i,n_j\}}$ belongs to part(${G\!\Sd_{B,A}}$) iff
either

\nav{}{$i$ is $s$ and $j$ is $t$, while ${m\in\{2GB,\ldots,2GB\pl
GA\mn 1\}}$, ${n<GA}$ and ${m\mn 2GB=n}$, or}

\vspace{-1ex}

\nav{}{$i$ and $j$ are both $s$, while ${m,n<2GB}$ and ${|m\mn
n|=GB}$.}

\noindent For $A$ and $B$ being as in the previous example, we
have

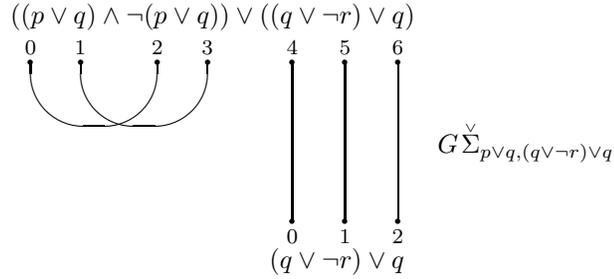
\begin{center}
\begin{picture}(180,110)
\put(100,20){\line(0,1){60}} \put(120,20){\line(0,1){60}}
\put(140,20){\line(0,1){60}}

\put(1,80){\circle*{2}} \put(20,80){\circle*{2}}
\put(49,80){\circle*{2}} \put(68,80){\circle*{2}}
\put(100,80){\circle*{2}} \put(120,80){\circle*{2}}
\put(140,80){\circle*{2}}

\put(100,20){\circle*{2}} \put(120,20){\circle*{2}}
\put(140,20){\circle*{2}}

\put(1,83){\makebox(0,0)[b]{\footnotesize 0}}
\put(20,83){\makebox(0,0)[b]{\footnotesize 1}}
\put(49,83){\makebox(0,0)[b]{\footnotesize 2}}
\put(68,83){\makebox(0,0)[b]{\footnotesize 3}}
\put(100,83){\makebox(0,0)[b]{\footnotesize 4}}
\put(120,83){\makebox(0,0)[b]{\footnotesize 5}}
\put(140,83){\makebox(0,0)[b]{\footnotesize 6}}

\put(100,17){\makebox(0,0)[t]{\footnotesize 0}}
\put(120,17){\makebox(0,0)[t]{\footnotesize 1}}
\put(140,17){\makebox(0,0)[t]{\footnotesize 2}}

\put(25,80){\oval(48,48)[b]} \put(44,80){\oval(48,48)[b]}

\put(155,50){\makebox(0,0)[l]{$G\!\Sd_{p\dis q,(q\dis\neg r)\dis
q}$}}

\put(117,10){\makebox(0,0)[t]{$(q\dis\neg r)\dis q$}}

\put(70,95){\makebox(0,0)[b]{$((p\dis q)\kon\neg(p\dis q))\dis
((q\dis\neg r)\dis q)$}}

\end{picture}
\end{center}

Let ${G(f\cirk g)=Gf\ast Gg}$. To define ${G(f\ks g)}$, for
${\!\ks\!\in\{\kon,\dis\}}$,  we need an auxiliary notion.

Suppose $b_X$ is a bijection from $X$ to $X_1$ and $b_Y$ a
bijection from $Y$ to $Y_1$. Then for ${R\subseteq(X^s\cup
Y^t)^2}$ we define ${R_{b_Y}^{b_X}\subseteq(X_1^s\cup Y_1^t)^2}$
by

\[(u_i,v_j)\in R_{b_Y}^{b_X}
\quad\mbox{\rm{iff}}\quad (i(b_U^{-1}(u)),j(b_V^{-1}(v)))\in R,
\]

\noindent where ${(i,U),(j,V)\in\{(s,X),(t,Y)\}}$.

If ${f\!:A\vdash D}$ and ${g\!:B\vdash E}$, then for
${\!\ks\!\in\{\kon,\dis\}}$ the set of ordered pairs ${G(f\ks g)}$
is

\[
Gf\cup Gg_{+GD}^{+GA}
\]

\noindent where ${+GA}$ is the bijection from $GB$ to $\{n\pl
GA\mid n\in GB\}$ that assigns ${n\pl GA}$ to $n$, and ${+GD}$ is
the bijection from $GE$ to $\{n\pl GD\mid n\in GE\}$ that assigns
${n\pl GD}$ to $n$.

It is not difficult to check that $G$ so defined is indeed a
functor from \PNN\ to \emph{Br}. For that, we determine by
induction on the length of derivation that for every equation
${f=g}$ of \PNN\ we have ${Gf=Gg}$ in \emph{Br}. We have shown by
this induction that \emph{Br} is a proof-net category, and the
existence of a structure-preserving functor $G$ from \PNN\ to
\emph{Br} follows from the freedom of \PNN.

We can define analogously to $G$ a functor, which we also call
$G$, from the category \mds\ to \emph{Br}. We just omit from the
definition of $G$ above the clauses involving $\Dk_{B,A}$ and
$\Sd_{B,A}$. The image of \mds\ by $G$ in \emph{Br} is the
subcategory of \emph{Br} isomorphic to \emph{Bij}, which we
mentioned above. The following is proved in \cite{DP04} (Section
7.6).

\prop{\mds\ Coherence}{The functor $G$ from \mds\ to Br is
faithful.}

\noindent It follows immediately from this coherence result that
\mds\ is isomorphic to a subcategory of \PNN\ (cf.\ \cite{DP04},
Section 14.4).

The following result is proved in \cite{DP05} (Section 2.7) and
\cite{DP05a}.

\prop{\PNN\ Coherence}{The functor $G$ from \PNN\ to Br is
faithful.}

\section{The category \SAps}

The objects of the category \SAps\ are the formulae of the
propositional language $\eL_{\top,\bot,\neg,\kon,\dis}$ generated
by \Pe, where $\neg$, $\kon$ and $\dis$ are as before, and $\top$
and $\bot$ are nullary connectives, i.e.\ propositional constants.
As primitive arrow terms we have $\mj_A$,
$\b{\kon}{\str}_{A,B,C}$, $\b{\kon}{\rts}_{A,B,C}$,
$\c{\kon}_{A,B}$, $\b{\dis}{\str}_{A,B,C}$,
$\b{\dis}{\rts}_{A,B,C}$, $\c{\dis}_{A,B}$, $d_{A,B,C}$ (see \S
2), $\Dk_{B,A}$, $\Sd_{B,A}$ (see \S 3), plus

\[
\begin{array}{ll}

\d{\kon}{\str}_A: A\kon\top \vdash A, & \d{\kon}{\rts}_A: A \vdash
  A\kon\top,
  \\[1ex]
  \d{\dis}{\str}_A: A\dis\bot \vdash A, & \d{\dis}{\rts}_A: A \vdash A\dis\bot,
 \\
\end{array}
\]

\noindent These primitive arrow terms together with the operations
on arrow terms $\cirk$, $\kon$ and $\dis$ (the same we had for
\mds\ and \PNN\ in \S\S2-3) define the arrow terms of \SAps.

The equations of \SAps\ are obtained by assuming all the equations
we have assumed for \PNN, plus

\begin{tabbing}
\quad\quad \= ($\d{\kon}{\str}$~{\it nat})\quad \= $f\cirk
\d{\kon}{\str}_A=\d{\kon}{\str}_B\cirk(f\kon\mj_\top)$,\kill

\> ($\d{\kon}{\str}$~{\it nat})\> $f\cirk\!
\d{\kon}{\str}_A\:=\:\d{\kon}{\str}_B\!\! \cirk(f\kon\mj_\top)$,
\\[1ex]
\> ($\d{\kon}{}\d{\kon}{}$)\>
$\d{\kon}{\str}_A\!\cirk\!\d{\kon}{\rts}_A\:=\mj_A$, \quad\quad
$\d{\kon}{\rts}_A\!\cirk\!\d{\kon}{\str}_A\:= \mj_{A\kon\top}$,
\\*[1ex]
\> ($\b{\kon}{}\d{\kon}{}$)\>
$\b{\kon}{\rts}_{A,B,\top}\!\cirk\!\d{\kon}{\rts}_{A\kon
B}\:=\mj_A\,\kon\d{\kon}{\rts}_B$,
\\[2ex]
\> ($\d{\dis}{\str}$~{\it nat})\> $f\cirk\!
\d{\dis}{\str}_A\;=\;\d{\dis}{\str}_B\!\cirk(f\dis\mj_\bot)$,
\\[1ex]
\> ($\d{\dis}{}\d{\dis}{}$)\>
$\d{\dis}{\str}_A\!\cirk\!\d{\dis}{\rts}_A\;=\mj_A$, \quad\quad
$\d{\dis}{\rts}_A\!\cirk\!\d{\dis}{\str}_A\;=\mj_{A\dis\bot}$,
\\*[1ex]
\> ($\b{\dis}{}\d{\dis}{}$)\>
$\b{\dis}{\rts}_{A,B,\bot}\!\cirk\!\d{\dis}{\rts}_{A\dis
B}\;=\mj_A\,\dis\d{\dis}{\rts}_B$,
\\[2ex]
for $\s{\kon}{\rts}_A\:=_{\df}\;\,
\c{\kon}_{A,\top}\!\cirk\!\d{\kon}{\rts}_A$,
\\[1ex]
\> ($d\!\s{\kon}{}$)\> $d_{\top,B,C}\,\cirk\!\s{\kon}{\rts}_{B\dis
C}\:$\= $=\;\s{\kon}{\rts}_B\dis\:\mj_C$,
\\*[1ex]
\> ($d\!\d{\dis}{}$)\> $\d{\dis}{\str}_{C\kon B}\!\cirk\,
d_{C,B,\bot}$\>$=\mj_C\:\kon\d{\dis}{\str}_B$,

\end{tabbing}

\noindent The set of equations of \SAps\ is closed under symmetry
and transitivity of equality and under the rules ({\it
cong}~$\ks\!$) for ${\!\ks\!\in\{\cirk,\kon,\dis\}}$ (see \S 2).
This defines the equations of \SAps.

We have the following definitions:

\[
\s{\kon}{\str}_A\:=_{\df}\;\,
\d{\kon}{\str}_A\!\cirk\!\c{\kon}_{\top,A}, \quad\quad
\s{\dis}{\str}_A\:=_{\df}\;\,
\d{\dis}{\str}_A\!\cirk\!\c{\dis}_{A,\bot}, \quad\quad
\s{\dis}{\rts}_A\:=_{\df}\;\,
\c{\dis}_{\bot,A}\!\cirk\!\d{\dis}{\rts}_A,
\]

\noindent which give isomorphisms in \SAps. Note that
${\s{\dis}{\str}_A:\bot\dis A\vdash A}$ is analogous to
$\Sd_{B,A}:(B\kon\neg B)\dis A\vdash A$, though $\Sd_{B,A}$ is not
an isomorphism. The equation (${\b{\kon}{}\Sd}$) of \S 3 is
analogous to the following equation of \SAps\ (an equation of
monoidal categories):

\[
\s{\dis}{\str}_{B\dis
A}\!\cirk\!\b{\dis}{\rts}_{\bot,B,A}\:=\;\s{\dis}{\str}_B\dis\:\mj_A.
\]

\noindent The equations (${d\!\s{\kon}{}}$) and
(${d\!\d{\dis}{}}$), which amount to the equations
(${\s{\kon}{}\!d^L}$) and (${\d{\dis}{}\!d^L}$) of Section 7.9 of
\cite{DP04} (these equations stem from \cite{CS97}, Section 2.1),
are analogous to the equations (${d\!\Sk}$) and (${d\!\Dd}$) of \S
3.

With the definitions

\begin{tabbing}

\quad\quad\quad\quad\quad\quad\quad\=$\tau^L_B\;\:$\=
$=_{\df}\;\;\,\s{\kon}{\str}_{\neg B\dis
B}\!\cirk\!\Dk_{B,\top}\,$\=$:\top\vdash \neg B\dis B$,
\\*[2ex]
\>$\gamma^R_B$\>$=_{\df}\;\;\,\Sd_{B,\bot}\!\cirk\!\d{\dis}{\rts}_{B\kon\neg
B}$\> $:B\kon \neg B\vdash\bot$,

\end{tabbing}

\noindent in \SAps, on the one hand, and

\begin{tabbing}

\quad\quad\quad\quad\quad\quad\quad\=$\Dk_{B,A}\;$\=
$=_{\df}\;\,(\mj_A\kon\tau^L_B)\cirk\!\d{\kon}{\rts}_A\:$\=$:A\vdash
A\kon(\neg B\dis B)$,
\\*[2ex]
\>$\Sd_{B,A}$\>$=_{\df}\;\;\,\s{\dis}{\str}_A\!\cirk(\gamma^R_B\dis\mj_A)$\>
$:(B\kon\neg B)\dis A\vdash A$,

\end{tabbing}

\noindent on the other hand, it can easily be established that
\SAps\ is isomorphic to the free \emph{symmetric linearly} (alias
\emph{weakly}) \emph{distributive category with negation} in the
sense of \cite{CS97} (Section 4, Definition 4.3) generated by \Pe.

\section{The Gentzenization of \SAps}

We will now define a new language of arrow terms to denote the
arrows of the category \SAps. We call these arrow terms {\it
Gentzen terms}, and we prove for Gentzen terms a result analogous
to Gentzen's cut-elimination theorem, which we will use to prove
that the category \PNN\ is isomorphic to a full subcategory of
\SAps.

As the arrow terms of \SAps, Gentzen terms will be defined
inductively starting from primitive Gentzen terms. As {\it
primitive} Gentzen terms we have $\mj_A\!:A\vdash A$, for $A$
being a letter, or $\top$, or $\bot$. To define the operations on
Gentzen terms, called {\it Gentzen operations}, which are mostly
partial operations, we need some preparation.

We define inductively a notion that for $\!\ks\!\in\{\kon,\dis\}$
we call a $\!\ks\!$-{\it context}:

\begin{itemize}
\item[] $\koc$ is a $\!\ks\!$-context; \item[] if $Z$ is a
$\!\ks\!$-context and $A$ an object of \SAps, then $Z\ks A$ and
$A\ks Z$ are $\!\ks\!$-contexts.
\end{itemize}

\noindent A $\!\ks\!$-context is called \emph{proper} when it is
not $\koc$.

Next we define inductively what it means for a $\!\ks\!$-context
$Z$ to be applied to an object $B$ of \SAps, which we write
$Z(B)$, or to an arrow term $f$ of \SAps, which we write $Z(f)$:

\begin{tabbing}
\mbox{\hspace{2em}}\= $(Z\ks A)(B)$ \= $=$ \=
\mbox{\hspace{10em}}\= $(Z\ks A)(f)$ \= $=$ \=\kill \> \>
$\koc(B)$\' $=$\> $B$,\> \> $\koc(f)$\' $=$\> $f$,
\\[1ex]
\> $(Z\ks A)(B)$\> $=$\> $Z(B)\ks A$,\> $(Z\ks A)(f)$\> $=$\>
$Z(f)\ks \mj_A$,
\\[1ex]
\> $(A\ks Z)(B)$\> $=$\> $A\ks Z(B)$;\> $(A\ks Z)(f)$\> $=$\>
$\mj_A\ks Z(f)$.
\end{tabbing}

\noindent We use $X$, perhaps with indices, as a variable for
$\kon$-contexts, and $Y$, perhaps with indices, as a variable for
$\dis$-contexts.

Then we have the Gentzen operation $\BE{\kon}{\rts}_X$, which
involves types specified by

\[
\f{f\!:X(A\kon(B\kon C))\vdash
D}{\BE{\kon}{\rts}_X\!\!f\!:X((A\kon B)\kon C)\vdash D}
\]

\noindent This is read ``if $f$ is a Gentzen term, then
$\BE{\kon}{\rts}_X\!\!f$ is a Gentzen term'', all that of the
required types. We use this rule notation for operations also in
the future. The Gentzen term $\BE{\kon}{\rts}_X\!\!f$ denotes the
arrow of \SAps\ named on the right-hand side of the $\den$ sign
below:

\[
\BE{\kon}{\rts}_X\!\!f\den f\cirk X(\b{\kon}{\rts}_{A,B,C}).
\]

We also have the following Gentzen operation:

\[
\f{f\!:D\vdash Y(A\dis(B\dis C))}{\BE{\dis}{\str}_Y\!\!f\den
Y(\b{\dis}{\str}_{A,B,C})\cirk f\!:D\vdash Y((A\dis B)\dis C)}
\]

\noindent and the following four analogous Gentzen operations,
where the types can be easily guessed:

\begin{tabbing}
\mbox{\hspace{10ex}}\= $\BE{\kon}{\str}_X\!\!f$ \= $\den$ \=
$f\cirk X(\b{\kon}{\str}_{A,B,C})$,\quad\quad\quad\=
$\BE{\dis}{\rts}_Y\!\!f$ \= $\den$ \=
$Y(\b{\dis}{\rts}_{A,B,C})\cirk f$,
\\[1ex]
\> $\CE{\kon}_X\!\!f$\> $\den$\> $f\cirk X(\c{\kon}_{A,B})$,\>
$\CE{\dis}_Y\!\!f$ \> $\den$ \> $Y(\c{\dis}_{A,B})\cirk f$.
\end{tabbing}

We also have the Gentzen operations in the following list:

\[
\f{f\!:A\vdash B}{\TE{\str}f\den f\cirk\s{\kon}{\str}_A:\top\kon
A\vdash B}\quad\quad\quad \f{f\!:B\vdash A}{\NE{\rts}f\den\;
\d{\dis}{\rts}_A\!\cirk f\!:B\vdash A\dis \bot}
\]

\[
\f{g\!:\top\kon A\vdash B}{\TE{\rts}g\den
g\cirk\s{\kon}{\rts}_A:A\vdash B}\quad\quad\quad \f{g\!:B\vdash
A\dis\bot}{\NE{\str}g\den\; \d{\dis}{\str}_A\!\cirk g\!:B\vdash A}
\]

\begin{tabbing}
for $\ep{\dis}_{D,C,B,A}=_{\df}\;(\c{\kon}_{C,D}\dis\:\mj_{B\dis
A})\cirk\!\b{\dis}{\rts}_{C\kon D,B,A}\!\cirk
((d_{C,D,B}\cirk\!\c{\kon}_{D\dis B,C})\dis \mj_A)\cirk$
\\*[.5ex]
\` $\cirk d_{D\dis B,C,A}:(D\dis B)\kon(C\dis A)\vdash(D\kon
C)\dis(B\dis A)$,
\end{tabbing}

\vspace{-2ex}

\[
\f{f_1\!:B_1\vdash A_1\dis C_1\quad\quad\quad\quad\quad\quad
f_2\!:B_2\vdash A_2\dis
C_2}{\Kon(f_1,f_2)\den\ep{\dis}_{A_1,A_2,C_1,C_2}\!\cirk(f_1\kon
f_2):B_1\kon B_2\vdash(A_1\kon A_2)\dis(C_1\dis C_2)}
\]

\begin{tabbing}
for $\ep{\kon}_{A,B,C,D}=_{\df}\; d_{A,C,B\kon
D}\cirk(\mj_A\kon(\c{\dis}_{C,B\kon D}\!\cirk
d_{B,D,C}))\cirk\!\b{\kon}{\rts}_{A,B,D\dis C}\!\cirk$
\\*[.5ex]
\` $\cirk(\mj_{A\kon B}\:\kon\c{\dis}_{D,C}):(A\kon B)\kon(C\dis
D)\vdash(A\kon C)\dis(B\kon D)$,
\end{tabbing}

\vspace{-2ex}

\[
\f{f_1\!:C_1\kon A_1\vdash B_1\quad\quad\quad\quad\quad\quad
f_2\!:C_2\kon A_2\vdash B_2}{\Dis(f_1,f_2)\den(f_1\dis
f_2)\cirk\ep{\kon}_{C_1,C_2,A_1,A_2}\!:(C_1\kon C_2)\kon(A_1\dis
A_2)\vdash B_1\dis B_2}
\]

\noindent (see \cite{DP04}, Section 7.6, for $\ep{\dis}$ and
$\ep{\kon}$),

\[
\f{f\!:B\vdash A\dis C}{\NL\!f\den\;\Sdp_{A,C}\!\cirk d_{\neg
A,A,C}\cirk\!\c{\kon}_{A\dis C,\neg A}\!\cirk(f\kon\mj_{\neg
A})\!:B\kon\neg A\vdash C}
\]

\[
\f{f\!:C\kon A\vdash B}{\NR\!f\den(\mj_{\neg A}\dis
f)\cirk\!\c{\dis}_{\neg A,C\kon A}\!\cirk d_{C,A,\neg
A}\cirk\!\Dkp_{A,C}\: :C\vdash\neg A\dis B}
\]

To define the remaining Gentzen operations, we need some
preparation. For every proper $\kon$-context $X$ we define
inductively as follows an object $E_X$ of \SAps:

\begin{tabbing}
\mbox{\hspace{2em}}\= $E_{\koci\kon B}$ \= = \=
$E_{B\kon\koci}=B$,\=
\\[1ex]
\> $E_{X\kon B}$\> =\> $E_X\kon B$,\> for $X$ proper,
\\[1ex]
\> $E_{B\kon X}$\> =\> $B\kon E_X$,\> for $X$ proper.
\end{tabbing}

\noindent For every proper $\kon$-context $X$ and every object $A$
of \SAps\ we define inductively as follows an arrow term
${\ta{\kon}_{X,A}:E_X\kon A\vdash X(A)}$ \SAps:

\begin{tabbing}
\mbox{\hspace{2em}}\= $\ta{\kon}_{B\kon\koci,A}$ \= $=_{\df}$ \=
$\mj_{B\kon A}:B\kon A\vdash B\kon A$,
\\[1.5ex]
\> $\ta{\kon}_{B\kon X,A}$\> $=_{\df}$\>
$(\mj_B\:\kon\ta{\kon}_{X,A})\cirk\!\b{\kon}{\rts}_{B,E_X,A}\:
:(B\kon E_X)\kon A\vdash B\kon X(A)$,
\\*[.5ex]
\` for $X$ proper,
\\[1ex]
\> $\ta{\kon}_{\koci\kon B,A}$\> $=_{\df}$\> $\c{\kon}_{B,A}\:
:B\kon A\vdash A\kon B$,
\\[1.5ex]
\> $\ta{\kon}_{X\kon B,A}$\> $=_{\df}$\>
$(\ta{\kon}_{X,A}\kon\:\mj_B)\cirk\!\b{\kon}{\str}_{E_X,A,B}\!
\cirk(\mj_{E_X}\:\kon\c{\kon}_{B,A})\cirk
\!\b{\kon}{\rts}_{E_X,B,A}:$
\\*[.5ex]
\` $(E_X\kon B)\kon A\vdash X(A)\kon B$, \quad for $X$ proper.
\end{tabbing}

For every proper $\dis$-context $Y$ we define inductively as
follows an object $D_Y$ of \SAps:

\begin{tabbing}
\mbox{\hspace{2em}}\= $D_{\koci\dis B}$ \= = \=
$D_{B\dis\koci}=B$,\=
\\[1ex]
\> $D_{Y\dis B}$\> =\> $D_Y\dis B$,\> for $Y$ proper,
\\[1ex]
\> $D_{B\dis Y}$\> =\> $B\dis D_Y$,\> for $Y$ proper.
\end{tabbing}

\noindent For every proper $\dis$-context $Y$ and every object $A$
of \SAps\ we define inductively as follows an arrow term
${\ta{\dis}_{Y,A}:Y(A)\vdash A\dis D_Y}$ of \SAps:

\begin{tabbing}
\mbox{\hspace{2em}}\= $\ta{\dis}_{\koci\dis B,A}$ \= $=_{\df}$ \=
$\mj_{A\dis B}:A\dis B\vdash A\dis B$,
\\[1.5ex]
\> $\ta{\dis}_{Y\dis B,A}$\> $=_{\df}$\>
$\b{\dis}{\rts}_{A,D_Y,B}\!\cirk(\ta{\dis}_{Y,A}\dis\:\mj_B):Y(A)\dis
B\vdash A\dis(D_Y\dis B)$,
\\*[.5ex]
\` for $Y$ proper,
\\[1ex]
\> $\ta{\dis}_{B\dis\koci,A}$\> $=_{\df}$\> $\c{\dis}_{A,B}\:
:B\dis A\vdash A\dis B$,
\\[1.5ex]
\> $\ta{\dis}_{B\dis Y,A}$\> $=_{\df}$\>
$\b{\dis}{\rts}_{A,B,D_Y}\!\!\!
\cirk(\c{\dis}_{A,B}\dis\:\mj_{D_Y})\cirk
\!\b{\dis}{\str}_{B,A,D_Y}\!\!\!\cirk(\mj_B\:\dis\ta{\dis}_{Y,A}):$
\\*[.5ex]
\` $B\dis Y(A)\vdash A\dis(B\dis D_Y)$,\quad for $Y$ proper.
\end{tabbing}

For ${f\!:A\vdash B}$, the following equations hold in \SAps:

\begin{tabbing}
\mbox{\hspace{4em}}\= $\mbox{($\ta{\kon}$ {\it nat})}$\quad\quad\=
$X(f)\cirk\!\ta{\kon}_{X,A}\;=\;\ta{\kon}_{X,B}\!\cirk(\mj_{E_X}\kon
f)$,
\\*[1.5ex]
\> $\mbox{($\ta{\dis}$ {\it nat})}$\>
$(f\dis\mj_{D_Y})\cirk\ta{\dis}_{Y,A}\;=\;\ta{\dis}_{Y,B}\!\cirk
Y(f)$;
\end{tabbing}

\noindent they are proved by applying naturality equations.

It is clear that for $\!\ks\!\in\{\kon,\dis\}$ and
$\ta{\xi}_{X,A}:A_1\vdash A_2$ there is an arrow term
$\tai{\xi}_{X,A}:A_2\vdash A_1$ of \SAps, which is a ``mirror
image'' of $\ta{\xi}_{X,A}$, such that in \SAps\ we have

\[
\tai{\xi}_{X,A}\!\cirk\!\ta{\xi}_{X,A}\;=\mj_{A_1},\quad\quad\quad
\ta{\xi}_{X,A}\!\cirk\!\tai{\xi}_{X,A}\;=\mj_{A_2}.
\]

\noindent For example, with

\begin{tabbing}
\mbox{\hspace{2em}}\= $\ta{\kon}_{F\kon((C\kon\koci)\kon B),A}$ \=
= \=
$(\mj_F\kon(\b{\kon}{\str}_{C,A,B}\!\cirk(\mj_C\:\kon\c{\kon}_{B,A})\cirk\!\b{\kon}{\rts}_{C,B,A}))\cirk\!\b{\kon}{\rts}_{F,C\kon
B,A}$
\\[1ex]
we have
\\[1ex]
\> $\tai{\kon}_{F\kon((C\kon\koci)\kon B),A}$\> =\>
$\b{\kon}{\str}_{F,C\kon B,A}\!\cirk
(\mj_F\kon(\b{\kon}{\str}_{C,B,A}\!\cirk
(\mj_C\:\kon\c{\kon}_{A,B})\cirk\!\b{\kon}{\rts}_{C,A,B})).$
\end{tabbing}

\noindent Officially, $\tai{\xi}_{X,A}$ is defined inductively as
$\ta{\xi}_{X,A}$, in a dual manner.

Next, we introduce the following abbreviation:

\begin{tabbing}
\centerline{$d_{X,A,Y}=_{\df}\;\;
\tai{\dis}_{Y,X(A)}\!\cirk(\ta{\kon}_{X,A}\dis\:\mj_{D_Y})\cirk
d_{E_X,A,D_Y}\cirk(\mj_{E_X}\kon\ta{\dis}_{Y,A})\cirk\tai{\kon}_{X,Y(A)}\:
:$}
\\*[.5ex]
\` $X(Y(A))\vdash Y(X(A))$.
\end{tabbing}

\noindent When $X$ or $Y$ is $\koc$, then we assume that
$d_{X,A,Y}$ stands for $\mj_{X(Y(A))}$, which is of type
${X(Y(A))\vdash Y(X(A))}$, i.e.\ ${Y(A)\vdash Y(A)}$ or
${X(A)\vdash X(A)}$.

We can finally define the remaining Gentzen operations, which are
all of the following form:

\[
\f{g\!:B\vdash Y(A)\quad\quad\quad\quad\quad\quad f\!:X(A)\vdash
C}{\cut_{X,Y}(f,g)\den Y(f)\cirk d_{X,A,Y}\cirk X(g):X(B)\vdash
Y(C)}
\]

\noindent This concludes the definition of Gentzen operations. The
set of Gentzen terms is the smallest set containing primitive
Gentzen terms and closed under the Gentzen operations above.

It is easy to infer from \mds\ Coherence of \S 4 that the
following equations hold in \SAps:

\begin{tabbing}
\mbox{\hspace{2em}}\= $\mbox{($d\Kon X$)}$\quad\quad\= $d_{A\kon
X,C,Y}$ \= = \= $d_{A\kon\koci,X(C),Y}\cirk(\mj_A\kon d_{X,C,Y})$,
\\[1ex]
\> $\mbox{($dX\Kon$)}$\> $d_{X\kon A,C,Y}$\> =\> $d_{\koci\kon
A,X(C),Y}\cirk(d_{X,C,Y}\kon \mj_A)$,
\\[1ex]
\> $\mbox{($d\Dis Y$)}$\> $d_{X,C,A\dis Y}$\> =\> $(\mj_A\dis
d_{X,C,Y})\cirk d_{X,Y(C),A\dis\koci}$,
\\[1ex]
\> $\mbox{($dY\Dis$)}$\> $d_{X,C,Y\dis A}$\> =\>
$(d_{X,C,Y}\dis\mj_A)\cirk d_{X,Y(C),\koci\dis A}$.
\end{tabbing}

\noindent The equation $\mbox{($d\Kon X$)}$ is analogous to the
equation $(d \kon)$ of \S 2, while $\mbox{($d\Dis Y$)}$ is
analogous to $(d \dis)$ of \S 2.

We can then prove the following

\prop{Gentzenization Lemma}{Every arrow of \SAps\ is denoted by a
Gentzen term.}

\dkz We first show by induction on the complexity of $A$ that for
every $A$ the arrow ${\mj_A\!:A\vdash A}$ is denoted by a Gentzen
term. For $A$ being a letter, or $\top$, or $\bot$, this is
trivial. For the induction step we use the following equations of
\SAps:

\begin{tabbing}
\mbox{\hspace{2em}}\= $\mbox{($\kon$)}$\quad\quad\=
$\:\:\NE{\str}\!\NE{\str}\!\BE{\dis}{\str}_{\koci}\!
\Kon(\NE{\rts}f_1,\NE{\rts}f_2)\,$ \= = \= $f_1\kon f_2$,
\\[1ex]
\> $\mbox{($\dis$)}$\>
$\TE{\rts}\TE{\rts}\!\BE{\kon}{\str}_{\koci}\!\Dis(\TE{\str}f_1,\TE{\str}f_2)$\>
=\> $f_1\dis f_2$.
\end{tabbing}

\noindent For $\mbox{($\kon$)}$ we use

\[
\ep{\dis}_{A_1,A_2,\bot,\bot}=(\mj_{A_1\kon
A_2}\:\dis\d{\dis}{\rts}_{\bot})\cirk\!\d{\dis}{\rts}_{A_1\kon
A_2}\! \cirk(\d{\dis}{\str}_{A_1}\kon\d{\dis}{\str}_{A_2}),
\]

\noindent which follows essentially from
(${\b{\dis}{}\d{\dis}{}}$) and (${d\!\d{\dis}{}}$) of \S 5 (we may
apply here the Symmetric Bimonoidal Coherence of \cite{DP04},
Section 6.4, which reduces to Mac Lane's symmetric monoidal
coherence of \cite{ML63}; see \cite{ML71}, Section VII.7, and
\cite{DP04}, Section 5.3). We proceed analogously for
$\mbox{($\dis$)}$.

We also have for the induction step the following equations of
\SAps:

\[
\NE{\str}\NR\CE{\kon}_{\koci}\!\NL\NE{\rts}\mj_A=
\TE{\rts}\NL\CE{\dis}_{\koci}\!\NR\TE{\str}\mj_A=\mj_{\neg A},
\]

\noindent for which we use (${d\!\d{\dis}{}}$) and
$\mbox{($\Sdp\Dkp$)}$, among other equations. The Gentzen term
that denotes $\mj_A$ is written $\mj_A$.

Next we have the following in \SAps:

\begin{tabbing}
\mbox{\hspace{2em}}\= $\BE{\kon}{\str}_{\koci}\!\!\mj_{(A\kon
B)\kon C}$ \= $\den$ \= $\b{\kon}{\str}_{A,B,C}$,\quad\quad\quad
\= $\BE{\dis}{\str}_{\koci}\!\!\mj_{A\dis(B\dis C)}$ \= $\den$ \=
$\b{\dis}{\str}_{A,B,C}$,
\\*[1ex]
\> $\BE{\kon}{\rts}_{\koci}\!\!\mj_{A\kon(B\kon C)}$\> $\den$\>
$\b{\kon}{\rts}_{A,B,C}$,\>
$\BE{\dis}{\rts}_{\koci}\!\!\mj_{(A\dis B)\dis C}$\> $\den$\>
$\b{\dis}{\rts}_{A,B,C}$,
\\[1.5ex]
\> \> $\CE{\kon}_{\koci}\!\!\mj_{B\kon A}$\' $\den$\>
$\c{\kon}_{A,B}$,\> \>$\CE{\dis}_{\koci}\!\!\mj_{B\dis A}$\'
$\den$\> $\c{\dis}_{A,B}$,
\\[1.5ex]
\pushtabs \> \hspace{4em} $\cut_{A\kon\koci,\koci\dis
C}(\mj_{A\kon B},\mj_{B\dis C})$ \= $\den$ \= $d_{A,B,C}$;
\\[2ex]
by using abbreviations according to $\mbox{($\kon$)}$ and
$\mbox{($\dis$)}$ above,
\\*[1.5ex]
\>
\hspace{4em}$\;\;\;\TE{\rts}\CE{\kon}_{\koci}\!(\mj_A\kon\NR\TE{\str}\mj_B)$\>
$\den$\> $\Dk_{B,A}$,
\\*[1ex]
\>
\hspace{4em}$\;\;\;\;\NE{\str}\CE{\dis}_{\koci}\!(\NL\NE{\rts}\mj_B\dis
\mj_A)$\> $\den$\> $\Sd_{B,A}$,
\\[1.5ex]
\poptabs \> \> $\CE{\kon}_{\koci}\!\TE{\str}\mj_A$\' $\!\den$\>
$\d{\kon}{\str}_A$,\> \>$\NE{\str}\mj_{A\dis \bot}$\' $\!\den$\>
$\d{\dis}{\str}_A$,
\\*[1ex]
\> \> $\TE{\rts}\CE{\kon}_{\koci}\!\mj_{A\kon \top}$\' $\!\den$\>
$\d{\kon}{\rts}_A$,\> \>$\NE{\rts}\mj_A$\' $\!\den$\>
$\d{\dis}{\rts}_A$.
\end{tabbing}

\noindent (For the equations involving $\Dk_{B,A}$ and $\Sd_{B,A}$
we rely on (${d\!\s{\kon}{}}$) and (${d\!\d{\dis}{}}$) of \S 5,
and on other equations, called stem-increasing equations in
\cite{DP05}, Section 2.5, and \cite{DP05a}, Section 6.)

For composition we have the following equation of \SAps:

\[
\cut_{\koci,\koci}(f,g)=f\cirk g,
\]

\noindent and for the operations $\kon$ and $\dis$ on arrows we
have the equations $\mbox{($\kon$)}$ and $\mbox{($\dis$)}$ above.
\qed

\vspace{2ex}

\section{Cut elimination in \SAps}

For the proof of the Cut-Elimination Theorem below we will
introduce analogues of Gentzen's notions of rank and degree. We
need some preliminary definitions to define these notions.

For $\!\ks\!\in\{\kon,\dis\}$, we define first by induction the
notion of $\!\ks\!$-{\it superficial} subformula of a formula of
$\eL_{\top,\bot,\neg,\kon,\dis}$:

\begin{itemize}
\item[] $A$ of the form $p$, $\bot$, $A_1\dis A_2$, or $\neg A'$,
is a $\kon$-superficial subformula of $A$; \item[] $A$ of the form
$p$, $\top$, $A_1\kon A_2$, or $\neg A'$, is a $\dis$-superficial
subformula of $A$; \item[] if $A$ is a $\!\ks\!$-superficial
subformula of $B$, then $A$ is a $\!\ks\!$-superficial subformula
of $B\ks C$ and $C\ks B$.
\end{itemize}

Consider a Gentzen term $f$ of the form

\[
\Kon(f_1,f_2)\!:B_1\kon B_2\vdash(A_1\kon A_2)\dis(C_1\dis C_2).
\]

\noindent The $\dis$-superficial subformula $A_1\kon A_2$ that is
the left disjunct of the target of $f$ is called the {\it leaf} of
$f$. All the other $\dis$-superficial subformulae of the target of
$f$, which are subformulae of $C_1$ or $C_2$, and all the
$\kon$-superficial subformulae of the source of $f$, which are
subformulae of $B_1$ or $B_2$, are called {\it lower parameters}
of $f$.

To every lower parameter $x$ of $f$, there corresponds
unambiguously a subformula $y$ in the target or the source of
either ${f_1\!:B_1\vdash A_1\dis C_1}$ or ${f_2\!:B_2\vdash
A_2\dis C_2}$, which we call the {\it upper parameter of f
corresponding to x}. The lower parameter $x$ is a
$\kon$-superficial subformula of the source of $f$ iff the
corresponding upper parameter $y$ is a $\kon$-superficial
subformula of the source of either $f_1$ or $f_2$ (it cannot be in
both), and analogously for parameters that are $\dis$-superficial
subformulae of targets. If $y$ is in the type of $f_1$, then $f_1$
is called the {\it subterm of f for the upper parameter y}, and
analogously for $f_2$.

For example, if $f$ is

\[
\Kon(\mj_{p\dis q},\NE{\rts}\mj_r)\!:(p\dis q)\kon r\vdash(p\kon
r)\dis(q\dis\bot),
\]

\noindent then $p\kon r$ in the target is the leaf of $f$, while
$q$ in the target of $f$ and ${p\dis q}$ and $r$ in the source of
$f$ are lower parameters of $f$. To the lower parameter $q$ of $f$
corresponds the upper parameter of $f$ that is the occurrence of
$q$ in the target of the subterm ${\mj_{p\dis q}\!:p\dis q\vdash
p\dis q}$ for this upper parameter; to the lower parameter ${p\dis
q}$ of $f$ corresponds the upper parameter of $f$ that is the
source of the subterm $\mj_{p\dis q}$ for this upper parameter;
and to the lower parameter $r$ of $f$ corresponds the upper
parameter of $f$ that is the source of the subterm
${\NE{\rts}\mj_r\!:r\vdash r\dis\bot}$ for this upper parameter.
Note that the subformula $\bot$ in the target of $f$ is not a
$\dis$-superficial subformula of this target, and hence is not a
lower parameter of $f$.

If the Gentzen term $f$ is of the form

\[
\Dis(f_1,f_2)\!:(C_1\kon C_2)\kon(A_1\dis A_2)\vdash B_1\dis B_2,
\]

\noindent then the $\kon$-superficial subformula ${A_1\dis A_2}$
that is the right conjunct of the source of $f$ is the leaf of
$f$, while all the other $\kon$-superficial subformulae of the
source of $f$ and the $\dis$-superficial subformulae of the target
of $f$ are the lower parameters of $f$. The upper parameters of
$f$ corresponding to these lower parameters, and the subterms of
$f$ for these upper parameters, are defined analogously to what we
had in the previous case.

The leaf of ${\NL f\!:B\kon\neg A\vdash C}$ is the
$\kon$-superficial subformula $\neg A$ that is the right conjunct
of its source, while the leaf of ${\NR f\!:C\vdash \neg A\dis B}$
is the $\dis$-superficial subformula $\neg A$ that is the left
disjunct of its target. In both cases, the remaining
$\kon$-superficial subformulae of the source or the remaining
$\dis$-superficial subformulae of the target are lower parameters,
to whom correspond, analogously to what we had before, upper
parameters in the source or target of the subterm $f$ for these
upper parameters.

If our Gentzen term is of the form

\[
\BE{\kon}{\rts}_X\!\!f, \BE{\kon}{\str}_X\!\!f,
\BE{\dis}{\str}_Y\!\!f, \BE{\dis}{\rts}_Y\!\!f, \CE{\kon}_X\!\!f,
\CE{\dis}_Y\!\!f, \TE{\str}\!f, \TE{\rts}\!f, \NE{\rts}\!f,
\NE{\str}\!f,\mbox{\rm{ or }} \cut_{X,Y}(f,g),
\]

\noindent then it has no leaves, and all the $\kon$-superficial
subformulae of its source and all the $\dis$-superficial
subformulae of its target are lower parameters, to which upper
parameters correspond in an obvious manner.

Finally, the Gentzen term ${\mj_p\!:p\vdash p}$ has two leaves,
which are its source $p$ and its target $p$. There are no
parameters of $\mj_p$, neither lower nor upper. The Gentzen term
${\mj_{\top}\!:\top\vdash\top}$ has as its leaf the target $\top$,
and no parameters (the source $\top$ of $\mj_{\top}$ is not a
$\kon$-superficial subformula of itself). The Gentzen term
${\mj_{\bot}\!:\bot\vdash\bot}$ has as its leaf the source $\bot$,
and no parameters (the target $\bot$ of $\mj_{\bot}$ is not a
$\dis$-superficial subformula of itself).

Let $x$ be a $\kon$-superficial subformula of the source of a
Gentzen term $f$ or a $\dis$-superficial subformula of the target
of $f$. Then the {\it cluster} of $x$ in $f$ is a sequence of
occurrences of formulae defined inductively as follows:

\begin{itemize}
\item[] if $x$ is a leaf of $f$, then the cluster of $x$ in $f$ is
$x$, \item[] if $x$ is not a leaf of $f$, then $x$ is a lower
parameter of $f$, and for $y_1$ being the upper parameter of $f$
corresponding to $x$, take the cluster $y_1\ldots y_n$, where
$n\geq 1$, of $y_1$ in the proper subterm $f'$ of $f$ that is the
subterm of $f$ for the upper parameter $y_1$ (the sequence
$y_1\ldots y_n$ is already defined, by the induction hypothesis);
the cluster of $x$ in $f$ is the sequence $xy_1\ldots y_n$.
\end{itemize}

All occurrences of formulae in a cluster are $\!\ks\!$-superficial
subformulae for $\!\ks\!$ being one of $\kon$ and $\dis$. If
$\!\ks\!$ is $\kon$, then the cluster is a {\it source} cluster,
and if $\!\ks\!$ is $\dis$, then it is a {\it target} cluster.

A {\it cut} is a Gentzen term of the form ${\cut_{X,Y}(f,g)}$. For
${g\!:B\vdash Y(A)}$ and ${f\!:X(A)\vdash C}$ let the formula $A$
be called the {\it cut formula} of the cut ${\cut_{X,Y}(f,g)}$.
Let $x$ be the displayed occurrence of $A$ in the source $X(A)$ of
$f$, and let $s$ be the length of the cluster of $x$ in $f$ (we
write $s$ because we have here a source cluster). Let $y$ be the
displayed occurrence of $A$ in the target $Y(A)$ of $g$, and let
$t$ be the length of the cluster of $y$ in $g$ (we write $t$
because we have here a target cluster).

Depending on the form of $A$, we define a number $r$, which we
call the {\it rank} of the cut ${\cut_{X,Y}(f,g)}$. If the cut
formula $A$ is of the form $p$ or $\neg A'$, then

\[
\begin{array}{ll}
r=\min(s,t)\mn 1, & \mbox{\rm if $A$ is $p$,}
\\[1ex]
r=s\pl t\mn 2, & \mbox{\rm if $A$ is $\neg A'$.}
\end{array}
\]

\noindent (As a matter of fact, when $A$ is $p$, we could
stipulate that $r$ is either ${s\pl t\mn 2}$, as when it is $\neg
A'$, or ${s\mn 1}$, or ${t\mn 1}$, but the computation of rank we
have introduced makes the cut-elimination procedure run faster,
and does not complicate the proof.)

If the cut formula $A$ is of the form $\top$ or $A_1\kon A_2$,
then ${r=t\mn 1}$. If, finally, the cut formula $A$ is of the form
$\bot$ or ${A_1\dis A_2}$, then ${r=s\mn 1}$.

We define the {\it degree} $d$ of a cut as the number of
occurrences of $\kon$, $\dis$ and $\neg$ in its cut formula. The
{\it complexity} of a cut is the ordered pair $(d,r)$, where $d$
is its degree and $r$ its rank. The complexities of cuts are
lexicographically ordered (i.e., ${(d_1,r_1)<(d_2,r_2)}$ iff
${d_1<d_2}$, or ${d_1=d_2}$ and ${r_1< r_2}$).

A Gentzen term is called {\it cut-free} when no subterm of it is a
cut. A cut ${\cut_{X,Y}(f,g)}$ is {\it topmost} when $f$ and $g$
are cut-free. (Since in the proof below, we compute the rank only
for topmost cuts, our definition of cluster can be shortened a
little bit by not considering the parameters of cuts; but this is
not a substantial shortening.)

We can then prove the following.

\prop{Cut-Elimination Theorem}{For every Gentzen term $h$ there is
a cut-free Gentzen term $h'$ such that ${h=h'}$ in \SAps.}

\dkz It suffices to prove the theorem when $h$ is a topmost cut.
We proceed by induction on the complexity ${(d,r)}$ of this
topmost cut.

Suppose ${r=0}$ and ${d=0}$. Then $h$ can be of one of the
following forms:

\begin{tabbing}
\mbox{\hspace{8em}}\= $\cut_{X,\koci}(f,\mj_A)$\quad\= for $A$
being $p$ or $\top$,
\\[1ex]
\> $\cut_{\koci,Y}(\mj_A,g)$\> for $A$ being $p$ or $\bot$,
\\[1ex]
and we have in \SAps
\\[1ex]
\> $\cut_{X,\koci}(f,\mj_A)=f$,
\\[1ex]
\> $\cut_{\koci,Y}(\mj_A,g)=g$.
\end{tabbing}

\noindent This settles the basis of the induction.

Suppose ${r=0}$ and ${d>0}$. Then the cut formula must be of the
form ${A_1\kon A_2}$ or ${A_1\dis A_2}$ or $\neg A'$. In the first
case, for ${f\!:X(A_1\kon A_2)\vdash D}$, ${g_1\!:B_1\vdash
A_1\dis C_1}$ and ${g_2\!:B_2\vdash A_2\dis C_2}$ we have the
equation

\[\cut_{X,\koci\dis(C_1\dis
C_2)}(f,\Kon(g_1,g_2))=\; \BE{\dis}{\rts}_{\koci}\!\!
\cut_{X'',\koci\dis C_2}(\cut_{X',\koci\dis C_1}(f,g_1),g_2)
\]

\noindent where $X'(C)$ is ${X(C\kon A_2)}$ and $X''(C)$ is
${X(B_1\kon C)}$. To prove this equation we apply naturality
equations and \mds\ Coherence of \S 4.

The complexity of the topmost cut $\cut_{X',\koci\dis C_1}(f,g_1)$
is ${(d',r')}$ with ${d'<d}$, and we can apply the induction
hypothesis to obtain a cut-free Gentzen term $f'$ equal to it in
\SAps. The complexity of the topmost cut $\cut_{X'',\koci\dis
C_2}(f',g_2)$ is ${(d'',r'')}$ with ${d''<d}$, and we can again
apply the induction hypothesis.

In case the cut formula is ${A_1\dis A_2}$, we have an analogous
equation, for which we use again \mds\ Coherence, and we reason
analogously, applying the induction hypothesis twice.

In case the cut formula is $\neg A'$, for ${f\!:D\kon A'\vdash E}$
and ${g\!:B\vdash A'\dis C}$ we have the equation

\[
\cut_{B\kon\koci,\koci\dis E}(\NL g,\NR
f)=\;\CE{\dis}_{\koci}\CE{\kon}_{\koci}\!\cut_{D\kon\koci,\koci\dis
C}(f,g),
\]

\noindent which holds by naturality equations and \PNN\ Coherence
of \S 4. Then we apply the induction hypothesis to the topmost cut
on the right-hand side, which has a smaller degree.

Suppose now ${r>0}$. If $r$ was computed as ${s\mn 1}$, or as
${s\pl t\mn 2}$, where ${s>1}$, then we may apply equations of
\SAps\ of the following form

\[
(\ast)\quad \cut_{X,Y}(\gamma
f',g)=\gamma_1\ldots\gamma_n\cut_{X',Y}(f',g)
\]

\noindent for $\gamma$, ${\gamma_1,\ldots,\gamma_n}$ unary Gentzen
operations. If ${(d,r)}$ is the complexity of the topmost cut
${\cut_{X,Y}(\gamma f',g)}$, then the complexity of the topmost
cut $\cut_{X',Y}(f',g)$ is ${(d,r-1)}$, and so we may apply to it
the induction hypothesis.

If $\gamma$ is a unary Gentzen operation different from
$\TE{\str}$, $\TE{\rts}$, $\NE{\rts}$ and $\NE{\str}$, then so are
$\gamma_1,\ldots,\gamma_n$, and to prove $(\ast)$ we apply
naturality equations and \PNN\ Coherence (sometimes \mds\
Coherence suffices, depending on $\gamma$). We have analogous
equations involving binary Gentzen operations, which are proved
analogously, relying on \mds\ Coherence (cf.\ \cite{DP04}, Section
11.2, Case (6), where on p.\ 251, in the second line
${\KR(f,\cut(g,h))}$ should be replaced by ${\KR(g,(f,h))}$, and
in the third line ${\cut(g,h)}$ should be replaced by
${\cut(f,h)}$).

If $\gamma$ in $(\ast)$ is $\TE{\str}$, then ${n=1}$ and
$\gamma_1$ is $\TE{\str}$. To prove $(\ast)$, we then apply
essentially the equation

\[
Y(\s{\kon}{\str}_{X(A)})\cirk d_{T\kon
X,A,Y}=d_{X,A,Y}\cirk\!\s{\kon}{\str}_{X(Y(A))},
\]

\noindent which we obtain with the help of $\mbox{($d\Kon X$)}$ of
the preceding section, (${d\!\s{\kon}{}}$) of \S 3.3, and
$\mbox{($\ta{\dis}$~{\it nat})}$ of the preceding section (as a
matter of fact, we may apply here the Symmetric Bimonoidal
Coherence of \cite{DP04}, Section 6.4). We proceed analogously if
$\gamma$ is $\TE{\rts}$.

If $\gamma$ in $(\ast)$ is $\NE{\rts}$ or $\NE{\str}$, then we
apply essentially Mac Lane's symmetric monoidal coherence of
\cite{ML63} (see also \cite{ML71}, Section VII.7, and \cite{DP04},
Section 5.3).

If $r$ was computed as ${t\mn 1}$, or as ${s\pl t\mn 2}$, where
${t>1}$, then we proceed in a dual manner. Instead of $(\ast)$, we
have equations of \SAps\ of the following form:

\[
\cut_{X,Y}(f,\gamma g')=\gamma_1\ldots\gamma_n\cut_{X,Y'}(f,g').
\]

\noindent This concludes the proof of the theorem. \qed

\section{\SAcons\ Coherence}

There is a functor $G$ from the category \SAps\ to \emph{Br},
which is defined as the functor $G$ from \PNN\ to \emph{Br} (see
\S 4) with the additional clauses that say that ${G\alpha}$ is an
identity arrow of \emph{Br} for $\alpha$ being $\d{\xi}{\str}_A$
and $\d{\xi}{\rts}_A$, where ${\!\ks\!\in\{\kon,\dis\}}$. It
follows from the existence of these functors and \PNN\ Coherence
of \S 4 that \PNN\ is isomorphic to a subcategory of \SAps\ (cf.\
\cite{DP04}, Section 14.4).

The following theorem can be proved with the help of the
Cut-Elimination Theorem of the preceding section.

\prop{Conservativeness Theorem}{If $A$ and $B$ are objects of
\PNN, then for every arrow ${f\!:A\vdash B}$ of \SAps\ there is an
arrow term ${f'\!:A\vdash B}$ of \PNN\ such that ${f=f'}$ in
\SAps.}

\noindent This theorem implies that \PNN\ is isomorphic to a full
subcategory of \SAps. In these isomorphisms every object of \PNN\
is mapped to itself, and so every object of \PNN\ in \SAps\ is in
the image of \PNN.

Let $\SA'$ be the full subcategory of \SA\ whose objects are all
the objects $A$ of \SA\ such that there is an isomorphism of type
${A\vdash A'}$ of \SA\ for $A'$ an object of \PNN. Then we can
restrict the functor $G$ from \SA\ to \emph{Br} to a functor $G$
from $\SA'$ to \emph{Br}, for which we can prove the following,
relying on the Conservativeness Theorem.

\prop{$\SA'$ Coherence}{The functor $G$ from $\SA'$ to Br is
faithful.}

\dkz Suppose $A$ and $B$ are objects of $\SA'$, and let
${j_A\!:A\vdash A'}$ and ${j_B\!:B\vdash B'}$ be isomorphisms of
\SA\ for $A'$ and $B'$ objects of \PNN. Suppose that
${f_1,f_2\!:A\vdash B}$ are arrows of \SA, i.e.\ of $\SA'$, such
that ${Gf_1=Gf_2}$.

Since \PNN\ is isomorphic to a full subcategory of \SA\ such that
every object of \PNN\ in \SA\ is in the image of \PNN, we have in
\SA\ that

\[
j_B\cirk f_i\cirk j^{-1}_A=f_i'
\]

\noindent for $i\in\{1,2\}$ and $f_i'$ an arrow term of \PNN. It
follows that ${Gf_1'=Gf_2'}$, and, according to what we said
immediately after the definition of the functor $G$ from \SA\ to
\emph{Br}, by \PNN\ Coherence we have that ${f_1'=f_2'}$ in \PNN,
and hence also in \SA. So ${f_1=f_2}$ in \SA. \qed

\vspace{2ex}

The category $\SA'$ is a category equivalent to \PNN, and its
coherence is a consequence of \PNN\ Coherence. We can find full
subcategories of $\SA'$ that are not only equivalent, but also
isomorphic, to \PNN.

Let \SAcons\ be the full subcategory of \SA\ whose objects are all
the objects $A$ of \SA\ such that there is an isomorphism of type
${A\vdash A'}$ of \SA\ for $A'$ being either an object of \PNN, or
$\top$, or $\bot$. Then we can restrict the functor $G$ from \SA\
to \emph{Br} to a functor $G$ from \SAcons\ to \emph{Br}, for
which we can prove the following, relying on the Conservativeness
Theorem and on $\SA'$ Coherence.

\prop{\SAcons\ Coherence}{The functor $G$ from \SAcons\ to Br is
faithful.}

\dkz There is no arrow of type ${\top\vdash\bot}$ in \SA.
(Otherwise, classical propositional logic would be inconsistent.)
There is also no arrow of type ${\bot\vdash\top}$ in \SA. If
${f\!:\bot\vdash\top}$ were such an arrow, then we would have in
\SA\ the arrow

\[
((\d{\kon}{\str}_p\!\cirk(\mj_p\kon f))\dis\:\mj_q)\cirk
d_{p,\bot,q}\cirk(\mj_p\:\kon\s{\dis}{\rts}_q)\!:p\kon q\vdash
p\dis q.
\]

\noindent Hence, by the Conservativeness Theorem, there would be
an arrow term $f'\!:p\kon q\vdash p\dis q$ of \PNN, and that such
an $f'$ does not exist can be shown by appealing to the
connectedness condition of proof nets (see \cite{DR89}).

Suppose $A$ and $B$ are objects of \SAcons; so $A$ and $B$ are
isomorphic in \SA\ to respectively $A'$ and $B'$ each of which is
either an object of \PNN, or $\top$, or $\bot$. Suppose that
${f_1,f_2\!:A\vdash B}$ are arrows of \SA, i.e.\ of \SAcons, such
that ${Gf_1=Gf_2}$.

As we have seen above, it is excluded that one of $A'$ and $B'$ is
$\top$ while the other is $\bot$. If $A'$ and $B'$ are objects of
\PNN, then we apply $\SA'$ Coherence.

Let $\SA_{+p}$ be \SA\ generated by ${\Pe\cup\{p\}}$ for a letter
$p$ foreign to \Pe, and hence also to $A$ and $B$. Let $\SA'_{+p}$
be the $\SA'$ subcategory of $\SA_{+p}$. In the remaining cases,
if either $A'$ or $B'$ is $\top$, then
${G(f_1\kon\mj_p)=G(f_2\kon\mj_p)}$. It is easy to see that
${f_1\kon\mj_p,f_2\kon\mj_p\!:A\kon p\vdash B\kon p}$ are arrows
of $\SA'_{+p}$, and so ${f_1\kon\mj_p=f_2\kon\mj_p}$ in $\SA_{+p}$
by $\SA'$ Coherence applied to $\SA'_{+p}$. Then in \SA\ generated
by \Pe\ we have ${f_1\kon\mj_{\top}=f_2\kon\mj_{\top}}$ (we just
substitute $\top$ for $p$ in the derivation of
$f_1\kon\mj_p=f_2\kon\mj_p$ in $\SA_{+p}$), and so we have in \SA\

\begin{tabbing}

\quad\quad\quad\quad\quad\quad\quad\quad\quad\quad\quad
$f_1\:$\=$=f_1\cirk\!\d{\kon}{\str}_A\!\cirk\!\d{\kon}{\rts}_A$,
\quad by $(\d{\kon}{}\d{\kon}{})$,
\\*[.5ex]
\>$=\;\d{\kon}{\str}_B\!\cirk(f_1\kon\mj_{\top})\cirk\!\d{\kon}{\rts}_A$,\quad
by ($\d{\kon}{\str}$~{\it nat}),
\\[.5ex]
\>$=\;\d{\kon}{\str}_B\!\cirk(f_2\kon\mj_{\top})\cirk\!\d{\kon}{\rts}_A$
\\[.5ex]
\>$=f_2$.

\end{tabbing}

\noindent If either $A'$ or $B'$ in the remaining cases is $\bot$,
then ${G(f_1\dis\mj_p)=G(f_2\dis\mj_p)}$, and we proceed
analogously. \qed

\vspace{2ex}

Let $\eL_{\top,\kon,\str}$ be the propositional language generated
by \Pe\ with the nullary connective $\top$ and the binary
connectives $\kon$ and $\str$. The formulae of
$\eL_{\top,\kon,\str}$ are the objects of the free symmetric
monoidal closed category \SMC\ generated by \Pe\ (see \cite{ML71},
Section VII.7, and \cite{DP05}, Section 3.1).

We call a formula $A$ of $\eL_{\top,\kon,\str}$
\emph{consequential} when for every subformula ${B\str C}$ of $A$
we have that either $B$ is letterless or $C$ has letters occurring
in it. An alternative way to characterize consequential formulae
is to say that these are formulae $A$ of $\eL_{\top,\kon,\str}$
for which there is an isomorphism of type ${A\vdash A'}$ of \SMC\
such that either $\top$ does not occur in $A'$ or $A'$ is $\top$.
(To establish the equivalence of these two characterizations, one
may rely on the results of \cite{DP97}.)

Let \SMCcons\ be the full subcategory of \SMC\ whose objects are
consequential formulae. With an appropriate definition of the
functor $G$ from \SMCcons\ to \emph{Br}, Kelly's and Mac Lane's
coherence theorem for symmetric monoidal closed categories of
\cite{KML71} amounts to the assertion that the functor $G$ from
\SMCcons\ to \emph{Br} is faithful. Both $\SA'$ Coherence and
\SAcons\ Coherence are analogous to this result of Kelly and Mac
Lane. For \SAcons\ Coherence the analogy is complete.

The proof of the Conservativeness Theorem is accomplished with the
help of a technical lemma, for whose formulation we introduce the
following terminology.

An object of \SAps, i.e.\ a formula of
$\eL_{\top,\bot,\neg,\kon,\dis}$, is {\it constant-free} when
neither $\top$ nor $\bot$ occurs in it. In other words, the
constant-free objects of \SAps\ are the objects of \PNN.

An object of \SAps\ is called {\it literate} when at least one
letter occurs in it; otherwise, it is \emph{letterless}. Every
constant-free formula is literate (but not conversely).

For ${\!\ks\!\!\in\{\kon,\dis\}}$, we define inductively when a
formula of $\eL_{\top,\bot,\neg,\kon,\dis}$ is $\!\ks\!$-{\it
nice}:

\vspace{1ex}

\nav{}{$\top$ is $\kon$-nice and $\bot$ is $\dis$-nice;}

\vspace{-1ex}

\nav{} {constant-free objects of \SAps\ are $\!\ks\!$-nice;}

\vspace{-1ex}

\nav{}{if $A$ and $B$ are $\!\ks\!$-nice, then $A\ks B$ is
$\!\ks\!$-nice.}

\vspace{1ex}

For a $\!\ks\!$-nice formula $A$ we define inductively an arrow
term ${\ro{\xi}{}_A:A\vdash A^r}$ of \SAps\ such that $A^r$ is
constant-free if $A$ is literate, $A^r$ is $\top$ if $A$ is
letterless and $\kon$-nice, and $A^r$ is $\bot$ if $A$ is
letterless and $\dis$-nice:

\begin{tabbing}
\mbox{\hspace{10em}}$\ro{\xi}{}_{A\kst B}\;$\=
$=\d{\xi}{\str}_A\!\cirk(\ro{\xi}{}_A\ks\ro{\xi}{}_B)$\kill

\quad\quad\quad\quad $\ro{\kon}{}_{\top}\;=\mj_{\top}$,
\mbox{\hspace{2.1em}} $\ro{\dis}{}_{\bot}\;$ \> $=\mj_{\bot}$,
\mbox{\hspace{2.1em}} $\ro{\xi}{}_A\;=\mj_A$,\` for $A$
constant-free,
\\[1ex]
\mbox{\hspace{10em}}$\ro{\xi}{}_{A\kst
B}$\>$=\;\ro{\xi}{}_A\ks\ro{\xi}{}_B$,\quad \` for $A$ and $B$
literate,
\\[1ex]
\mbox{\hspace{9.7em}} $\ro{\xi}{}_{A\kst B}\;$\>
$=\;\d{\xi}{\str}_A\!\cirk(\ro{\xi}{}_A\ks\ro{\xi}{}_B)$,\` for
$B$ letterless,
\\*[1ex]
\mbox{\hspace{9.7em}} $\ro{\xi}{}_{A\kst B}\;$\>
$=\;\s{\xi}{\str}_B\!\cirk(\ro{\xi}{}_A\ks\ro{\xi}{}_B)$, \` for
$A$ letterless.
\end{tabbing}

\noindent It is clear that $\ro{\xi}{}_A$ is an isomorphism of
\SAps, with inverse ${\ro{\xi}{-1}_A:A^r\vdash A}$.

The Conservativeness Theorem is a corollary of the following lemma
(we just instantiate statement (1) of this lemma).

\prop{Lemma}{Let ${f\!:A\vdash B}$ be an arrow of \SAps\ such that
$A$ is $\kon$-nice and $B$ is $\dis$-nice.}

\vspace{-1ex}

\nav{(1)}{\it If both $A$ and $B$ are literate, then there is an
arrow term $f^r\!:A^r\vdash B^r$ of $\;$\PNN\ such that in \SAps\
we have}

\vspace{-2ex}

\[
\ro{\dis}{}_B\!\cirk f\cirk\!\ro{\kon}{-1}_A\;=f^r\!.
\]

\vspace{-1ex}

\nav{(2)}{\it If $A$ is letterless and $B$ is literate, then for
every constant-free $C$ there is an arrow term ${f^r\!:C\vdash
C\kon B^r}$ of $\;$\PNN\ such that in \SAps\ we have}

\vspace{-2ex}

\[
(\mj_C\kon(\ro{\dis}{}_B\!\cirk
f\cirk\!\ro{\kon}{-1}_A))\cirk\!\d{\kon}{\rts}_C\;=f^r\!.
\]

\vspace{-1ex}

\nav{(3)}{\it If $A$ is literate and $B$ is letterless, then for
every constant-free $C$ there is an arrow term ${f^r\!:A^r\dis
C\vdash C}$ of $\;$\PNN\ such that in \SAps\ we have}

\vspace{-2ex}

\[
\s{\dis}{\str}_C\!\cirk((\ro{\dis}{}_B\!\cirk f\cirk\!
\ro{\kon}{-1}_A)\dis\mj_C)=f^r\!.
\]

The proof of this lemma, which may be found in \cite{DP05}
(Section 4.3), is based on the Gentzenization Lemma and the
Cut-Elimination Theorem of the preceding two sections. We take
that $f$ in the lemma is a cut-free Gentzen term, and we proceed
by induction on the complexity of $f$.

\baselineskip=0.84\baselineskip

\end{document}